\documentclass[10pt]{article}
\usepackage{makeidx}
\usepackage{latexsym}
\usepackage{amsfonts}
\usepackage{amsmath}
\usepackage{amsthm}
\usepackage{amssymb}
\usepackage{amscd}
\usepackage[dvips]{graphicx}
\usepackage{color}
\usepackage{eepic}
\usepackage{setspace}
\usepackage{calrsfs}

\def\C{\mathbb{ C}}

\newdimen\Squaresize \Squaresize=20pt
\newdimen\thickness \thickness=1pt

\def\Square#1{\hbox{\vrule width \thickness
   \vbox to \Squaresize{\hrule height \thickness\vss
      \hbox to \Squaresize{\hss#1\hss}
   \vss\hrule height\thickness}
\unskip\vrule width \thickness}
\kern-\thickness}

\def\vsquare#1{\vbox{\Square{$#1$}}\kern-\thickness}
\def\blank{\omit\hskip\Squaresize}

\def\fibyoung#1{\let\\=\cr              
\vbox{\smallskip\offinterlineskip
\halign{&\vsquare{##}\cr #1}}\,}

\def\borderlessrect#1#2{\hbox{\hskip \thickness
   \vbox to \Squaresize{\vskip \thickness \vss
      \hbox to #2 {\hss #1\hss}
   \vss\vskip\thickness}
\unskip\hskip \thickness}
\kern-\thickness}

\def\vborderlessrect#1#2{\vbox{\borderlessrect{$#1$}{#2}}\kern-\thickness}

\def\borderless#1{\omit\vborderlessrect{#1}{\Squaresize}}
\def\borderlessrc#1#2{\omit\vborderlessrect{#1}{#2}}

\def\msquare#1{\vbox{\hbox{\vrule width \thickness
   \vbox to \Squaresize{\hrule height \thickness
      \hbox to \Squaresize{\hfil{\sevenrm #1}}
   \vfil\hrule height\thickness}
\unskip\vrule width \thickness}
\kern-\thickness}\kern-\thickness}

\def\twosquare#1#2{\vbox{\hbox{\vrule width \thickness 
   \vbox to \Squaresize{\hrule height \thickness
      \hbox to \Squaresize{\hfil{\rm #1}}\vss
      \hbox to \Squaresize{\hss{#2}\hss}
   \vfil\hrule height\thickness}
\unskip\vrule width \thickness}
\kern-\thickness}\kern-\thickness}

\def\twoblank#1#2{\vbox{\hbox{
   \vbox to \Squaresize{\vskip 2pt
      \hbox to \Squaresize{\hfil{\sevenrm #1}\ }\vss
      \hbox to \Squaresize{\hss{#2}\hss}
   \vfil}\unskip\kern-\thickness}
}\unskip\kern-\thickness}

\def\young#1{
\def\>{\blank}
\def\<{\borderless}
\def\*{\borderlessrc}
\def\p{\omit\msquare}
\def\t{\omit\twosquare}
\def\b{\omit\twoblank}
\let\\=\cr 
\vbox{\smallskip\offinterlineskip
\halign{&\vsquare{##}\cr #1}}}

\newdimen\smsquaresize \smsquaresize=12pt
\newdimen\smthickness \smthickness=.5pt
\font\smcellfont=cmss8 scaled \magstep0

\def\smsquare#1{\hbox{\vrule width \smthickness
   \unskip\vbox to \smsquaresize{\hrule height \smthickness\vss
      \hbox to \smsquaresize{\hss{\smcellfont #1}\hss}
   \vss\hrule height\smthickness}
\unskip\vrule width \smthickness}
\kern-\smthickness}

\def\smvsquare#1{\vbox{\smsquare{$#1$}}\kern-\smthickness}
\def\blank{\omit\hskip\smsquaresize}

\def\smyoung#1{\let\\=\cr 
\vbox{\smallskip\offinterlineskip
\halign{&\smvsquare{##}\cr #1}}}
\newdimen\vsmsquaresize \vsmsquaresize=10pt
\newdimen\vsmthickness \vsmthickness=.5pt
\font\vsmcellfont=cmsl8 scaled \magstep0
\font\vsmletterfont=cmr6 scaled \magstep0

\def\vsmsquare#1{\hbox{\vrule width \vsmthickness
   \unskip\vbox to \vsmsquaresize{\hrule height \vsmthickness\vss
      \hbox to \vsmsquaresize{\hss{\vsmcellfont #1}\hss}
   \vss\hrule height\vsmthickness}
\unskip\vrule width \vsmthickness}
\kern-\vsmthickness}
\def\vsmvsquare#1{\vbox{\vsmsquare{#1}}\kern-\vsmthickness}
\def\vsmblank{\omit\hskip\vsmsquaresize}
\def\vsmborderless#1{\hbox{\hskip \vsmthickness\unskip
   \vbox to \vsmsquaresize{\vss
      \hbox to \vsmsquaresize{\hss{\vsmletterfont #1}\hss}
   \vss}
\unskip\hskip \vsmthickness}
\kern-\vsmthickness}                                                            \def\vsmvborderless#1{\vbox{\vsmborderless{#1}}\kern-\vsmthickness}

\def\vsmyoung#1{
\def\>{\vsmblank}
\def\<{\omit\vsmvborderless}
\let\\=\cr 
\vbox{\smallskip\offinterlineskip
\halign{&\vsmvsquare{##}\cr #1}}}

\newdimen\edgesize \edgesize=20pt
\newdimen\eedgesize \eedgesize=21pt
\newdimen\doublesize \doublesize=41pt
\newdimen\ddoublesize \ddoublesize=43pt

\newdimen\triplesize \triplesize=62pt
\newdimen\tetrasize \tetrasize=85pt
\newdimen\futosa \futosa=1pt

\def\Hako#1{\hbox{\vrule width \futosa
   \vbox to \eedgesize{\hrule height \futosa\vss
      \hbox to \edgesize{\hss#1\hss}
   \vss\hrule height\futosa}
\unskip\vrule width \futosa}
\kern-\futosa}

\def\Seihokei#1{\vbox{\Hako{#1}}\kern-\futosa}

\def\Horizontal#1{\hbox{\vrule width \futosa
   \vbox to \eedgesize{\hrule height \futosa\vss
      \hbox to \doublesize{\hss#1\hss}
   \vss\hrule height \futosa}
\unskip\vrule width \futosa}
\kern-\futosa}

\def\Vertical#1{\hbox{\vrule width \futosa
   \vbox to \ddoublesize{\hrule height \futosa\vss
      \hbox to \edgesize{\hss#1\hss}
   \vss\hrule height \futosa}
\unskip\vrule width \futosa}
\kern-\futosa}

\def\NS#1{\vbox{\Hako{$#1$}\vskip22pt}\kern-\futosa}
\def\NSS#1{\vbox{\Hako{$#1$}\vskip42pt}\kern-\futosa}
\def\NSSS#1{\vbox{\Hako{$#1$}\vskip64pt}\kern-\futosa}
\def\H#1{\vbox{\Horizontal{$#1$}}\kern-\futosa}
\def\HH#1#2{\vbox{\Horizontal{$#1$}\Horizontal{$#2$}}\kern-\futosa}
\def\VER#1{\vbox{\Vertical{$#1$}}\kern-\futosa}
\def\VS#1{\vbox{\Vertical{$#1$}\vskip20pt}\kern-\futosa}
\def\VSS#1{\vbox{\Vertical{$#1$}\vskip42pt}\kern-\futosa}
\def\VV#1#2{\vbox{\Vertical{$#1$}\Vertical{$#2$}}\kern-\futosa}
\def\NV#1#2{\vbox{\Hako{$#1$}\Vertical{$#2$}}\kern-\futosa}
\def\NVS#1#2{\vbox{\Hako{$#1$}\Vertical{$#2$}\vskip22pt}\kern-\futosa}
\def\YTT#1#2#3{\vbox{\Horizontal{$#1$}\hbox{\vbox{\Vertical{$#2$}}\kern-\futosa\vbox{\Vertical{$#3$}}\kern-\futosa}}\kern-\futosa}

\def\domino#1{
\def\ns{\omit\NS}
\def\nss{\omit\NSS}
\def\nsss{\omit\NSSS}
\def\h{\omit\H}
\def\hh{\omit\HH}
\def\V{\omit\VER}
\def\Vs{\omit\VS}
\def\Vss{\omit\VSS}
\def\vsss{\omit\VSSS}
\def\vv{\omit\VV}
\def\nv{\omit\NV}
\def\nvs{\omit\NVS}
\def\ytt{\omit\YTT}
\let\\=\cr 
\vbox{\smallskip\offinterlineskip
\halign{&\Seihokei{##}\cr #1}}}

\theoremstyle{definition}
\newtheorem{theorem}{Theorem}[section]
\newtheorem{prop}[theorem]{Proposition}
\newtheorem{lemma}[theorem]{Lemma}
\newtheorem{corollary}[theorem]{Corollary}

\newtheorem{remark}[theorem]{Remark}

\newenvironment{demo}[1]{%
  \trivlist
  \item[\hskip\labelsep
        {\bf #1.}]
}{%
\hfill\qedsymbol
  \endtrivlist
}

\renewcommand{\mathcal}{\mathrsfs}

\newcommand{\p}{\partial}

\renewcommand{\b}{\beta}

\newcommand\sgn{\operatorname{sgn}}

\newcommand\rdots{\mathinner{\mkern1mu\raise0pt\vbox{\kern7pt\hbox{.}}
     \mkern2mu\raise4pt\hbox{.}\mkern2mu\raise8pt\hbox{.}\mkern1mu}}

\def\PATH#1#2{{\cal P}\left({#1},{#2}\right)}
\def\NPATH#1#2{{\cal P}_0\left({#1},{#2}\right)}
\def\GF#1{{\operatorname{GF}}\left[{#1}\right]}

\def\MOD{\operatorname{mod}}

\def\V#1#2#3#4{V^{{#1},{#2}}\left({#3};{#4}\right)}

\numberwithin{equation}{section}
\def\rdots{\mathinner{\mkern1mu\raise0pt\vbox{\kern7pt\hbox{.}}
     \mkern2mu\raise4pt\hbox{.}\mkern2mu\raise8pt\hbox{.}\mkern1mu}}

\def\covered{\mathinner{\mkern1mu\raise0pt\vbox{\kern7pt\hbox{$<$}}
     \mkern-4mu\raise2pt\hbox{.}\mkern2mu}}
\def\covers{\mathinner{\mkern1mu\raise0pt\vbox{\kern7pt\hbox{$>$}}
     \mkern-12mu\raise2pt\hbox{.}\mkern8mu}}

\def\defterm#1{{\sl #1}\/}
\def\newterm#1{{\sl #1}\/}

\def\operatorname#1{{\mathrm{#1}\>\!}}
\def\Cal#1{{\mathcal{#1}}}
\def\Bbb#1{{\mathbb{#1}}}
\def\qbinom#1#2{\left[{{#1}\atop{#2}}\right]_{q}}

\title{
A $q$-analogue of Catalan Hankel determinants
}

\author{
Masao Ishikawa\\
\small Faculty of Education, Tottori University\\[-0.8ex]
\small Koyama, Tottori, Japan\\[-0.8ex]
\small \texttt{ishikawa@fed.tottori-u.ac.jp}
\and
Hiroyuki Tagawa\\
\small Faculty of Education, Wakayama University\\[-0.8ex]
\small Sakaedani, Wakayama, Japan\\[-0.8ex]
\small \texttt{tagawa@math.edu.wakayama-u.ac.jp}
\and
Jiang Zeng\\
\small Institut Camille Jordan\\[-0.8ex]
\small Universit\'e Claude Bernard Lyon I\\[-0.8ex]
\small 43, boulevard du 11 novembre 1918\\[-0.8ex]
\small 69622 Villeurbanne Cedex, France\\[-0.8ex]
\small \texttt{zeng@math.univ-lyon1.fr}
}

\date{
}

\begin{document}

\maketitle

\abstract{
In this article
 we shall survey the various methods of evaluating Hankel determinants 
and as an illustration we evaluate some Hankel determinants of a $q$-analogue of Catalan numbers. 
Here we consider $\frac{(aq;q)_{n}}{(abq^{2};q)_{n}}$ as a $q$-analogue 
of Catalan numbers $C_{n}=\frac1{n+1}\binom{2n}{n}$,
which is known as the moments of the little $q$-Jacobi polynomials.
We also give several proofs of this $q$-analogue,
in which we use lattice paths, the orthogonal polynomials,
or the basic hypergeometric series.
We also consider a $q$-analogue of Schr\"oder Hankel determinants,
and give a new proof of Moztkin Hankel determinants using an addition formula for
${}_2F_{1}$.
}
\medskip

\section{
Introduction
}\label{sec:intro}


Given a sequence $a_{0}$, $a_{1}$, $a_{2}$,$\dots$,
we set the Hankel matrix of the sequence to be
\begin{equation}
A_{n}^{(t)}
=\left(a_{i+j+t}\right)_{0\leq i,j\leq n-1}
=\begin{pmatrix}
a_{t}&a_{t+1}&\hdots&a_{t+n-1}\\
a_{t+1}&a_{t+2}&\hdots&a_{t+n}\\
\vdots&\vdots&\ddots&\vdots\\
a_{t+n-1}&a_{t+n}&\hdots&a_{t+2n-2}\\
\end{pmatrix}.
\label{eq:Hankel}
\end{equation}
For $n=0,1,2,\dots$,
let
\begin{equation}
C_{n}=\frac1{n+1}\binom{2n}{n},
\label{eq:Catalan}
\end{equation}
which are called the \newterm{Catalan numbers}.
The generating function for the Catalan numbers is given by
\begin{equation*}
\sum_{n\geq0}C_{n}t^{n}
=\frac{1-\sqrt{1-4t}}{2t}.
\end{equation*}
If we put $a_{n}=C_{n}$ in \eqref{eq:Hankel},
then the following identity is well-known and several proofs are known \cite{BCQY,CRI,MW,T,V}:
\begin{equation}
\det A_{n}^{(t)}
=\det\left(C_{i+j+t}\right)_{0\leq i,j\leq n-1}
=\prod_{1\leq i\leq j\leq t-1}\frac{i+j+2n}{i+j}.
\label{eq:Catalan-det}
\end{equation}
If we put $B_{n}=\binom{2n+1}{n}$ and $D_{n}=\binom{2n}{n}$,
then the following variations are also known \cite{T}:
\begin{align}
&\det\left(B_{i+j+t}\right)_{0\leq i,j\leq n-1}
=\prod_{1\leq i\leq j\leq t-1}\frac{i+j-1+2n}{i+j-1},
\label{eq:B-Catalan}
\\
&\det\left(D_{i+j+t}\right)_{0\leq i,j\leq n-1}
=2^{n}\prod_{1\leq i< j\leq t-1}\frac{i+j+2n}{i+j}.
\label{eq:D-Catalan}
\end{align}
As a generalization of \eqref{eq:Catalan-det},
Krattenthaler \cite{K} has obtained
\begin{equation}
  \det\left(C_{k_{i}+j}\right)_{0\leq i,j\leq n-1}
  =\prod_{0\leq i<j\leq n-1}\left(k_{j}-k_{i}\right)
  \prod_{i=0}^{n-1}\frac{(i+n)!(2k_{i})!}{(2i)!k_{i}!(k_{i}+n)!}
\label{eq:Krattenthaler}
\end{equation}
for a positive integer $n$ and non-negative integers
$k_{0}$, $k_{1}$,$\dots$, $k_{n-1}$.

In this article we shall survey the various methods of evaluating Hankel determinants and as an illustration we give a q-analogue of the above results. We first recall some terminology in q-series
(see Gasper-Rahman's book \cite{GR}) before stating the main theorem.
 Next some terminology is defined before stating the main theorem.
We use the notation:
\begin{equation*}
    (a;q)_{\infty}=\prod_{k=0}^{\infty}(1-aq^{k}),\qquad
    (a;q)_{n}=\prod_{k=0}^{n-1}(1-aq^{k})
\end{equation*}
for a nonnegative integer $n\geq0$.
Usually  $(a;q)_{n}$ is called the  \defterm{$q$-shifted factorial},
and we frequently use the compact notation:
\begin{align*}
    &(a_{1},a_{2},\dots,a_{r};q)_{\infty}=(a_{1};q)_{\infty}(a_{2};q)_{\infty}\cdots(a_{r};q)_{\infty},\\
    &(a_{1},a_{2},\dots,a_{r};q)_{n}=(a_{1};q)_{n}(a_{2};q)_{n}\cdots(a_{r};q)_{n}.
\end{align*}
If we put $a=q^{\alpha}$ and $q\rightarrow1$, then we have
\begin{equation*}
\lim_{q\rightarrow1}\frac{(q^{\alpha};q)_{n}}{(1-q)^n}
=(\alpha)_{n},
\end{equation*}
where
$\displaystyle(\alpha)_{n}=\prod_{k=0}^{n-1}(\alpha+k)$
is called the \defterm{raising factorial}.
We shall define the \defterm{${}_{r+1}\phi_{r} $ basic hypergeometric series} by
\begin{align*}
{}_{r+1}\phi_{r}\left[\,
{{a_{1},a_{2},\dots,a_{r+1}}\atop{b_{1},\dots,b_{r}}};q,z
\,\right]
=\sum_{n=0}^{\infty}\frac{(a_{1},a_{2},\dots,a_{r+1};q)_{n}}{(q,b_{1},\dots,b_{r};q)_{n}}z^{n}.
\end{align*}
If we put $a_{i}=q^{\alpha_{i}}$ and $b_{i}=q^{\beta_{i}}$ in the above series
and let $q\rightarrow1$,
then we obtain the \defterm{${}_{r+1}F_{r} $ hypergeometric series}
\begin{align*}
{}_{r+1}F_{r}\left[\,
{{\alpha_{1},\alpha_{2},\dots,\alpha_{r+1}}\atop{\beta_{1},\dots,\beta_{r}}};z
\,\right]
=\sum_{n=0}^{\infty}\frac{(\alpha_{1})_{n}(\alpha_{2})_{n}\cdots(\alpha_{r+1})_{n}}{n!(\beta_{1})_{n}\dots(\beta_{r})_{n}}z^{n}.
\end{align*}
The Motzkin number $M_{n}$ is defined to be
\[
M_{n}={}_2F_{1}\left[{{(1-n)/2,-n/2}\atop{2}};4\right].
\]
The generating function for the Motzkin numbers is given by
\begin{align*}
\sum_{n=0}^{\infty}M_{n}x^n
=\frac{1-x-\sqrt{1-2x-3x^2}}{2x^2}.
\end{align*}
It is known \cite{A} that
\begin{equation}
\det\left(M_{i+j}\right)_{0\leq i,j\leq n-1}=1
\label{eq:Aigner01}
\end{equation}
for $n\geq1$,
and
\begin{equation}
\det \left(M_{i+j+1}\right)_{0\leq i,j\leq n-1}=1,0,-1
\label{eq:Aigner02}
\end{equation}
for $n\equiv0,1$ ($\operatorname{mod}6$),
$n\equiv2,5$ ($\operatorname{mod}6$),
$n\equiv3,4$ ($\operatorname{mod}6$),
respectively.

\smallskip
The large Schr\"oder number $S_{n}$ is defined to be
\[
S_{n}=2{}_2F_{1}\left[{{-n+1,n+2}\atop{2}};-1\right]
\]
for $n\geq1$ ($S_0=1$).
The generating function for the large Schr\"oder numbers is
\begin{align}
\sum_{n=0}^{\infty}S_{n}x^n
=\frac{1-x-\sqrt{1-6x+x^2}}{2x}.
\end{align}
Eu and Fu \cite{EF} have proved
\begin{equation}
    \det\left(S_{i+j}\right)_{0\leq i,j\leq n-1}=2^{\binom{n}2},
\qquad
    \det\left(S_{i+j+1}\right)_{0\leq i,j\leq n-1}=2^{\binom{n+1}2}
\label{eq:Schroder,t=0,1}
\end{equation}
  for $n\geq1$ (see \cite{BK,EF,SX}).
We can also prove that
\begin{equation}
    \det\left(S_{i+j+2}\right)_{0\leq i,j\leq n-1}=2^{\binom{n+1}2}(2^{n+1}-1)
\label{eq:Schroder,t=2}
\end{equation}
holds  for $n\geq1$.

\medbreak

In this article,
as a generalization of \eqref{eq:Catalan},
we choose
\begin{equation}
\mu_{n}=\frac{(aq;q)_{n}}{(abq^2;q)_{n}}
\label{eq:moment}
\end{equation}
for a nonnegative integer $n$.
The aim of this article is to give three different proofs of the following theorem:

\begin{theorem}\label{th:q-Catalan-Hankel}
Let $n$ be a positive integer.
Then we have
\begin{equation}
  \det\left(\mu_{i+j}\right)_{0\leq i,j\leq n-1}
  =a^{\frac12n(n-1)}q^{\frac16n(n-1)(2n-1)}
  \prod_{k=1}^{n}\frac{(q,aq,bq;q)_{n-k}}
  {(abq^{n-k+1};q)_{n-k}(abq^{2};q)_{2(n-k)}}.
\label{eq:q-catalan-det}
\end{equation}
\end{theorem}
As a corollary of this theorem we can get the following more general identity.
\begin{corollary}\label{cor:q-Catalan-Hankel}
Let $n$ be a positive integer,
and $t$ a nonnegative integer.
Then we have
  \begin{align}
  \det\left(\mu_{i+j+t}\right)_{0\leq i,j\leq n-1}
  &=a^{\frac12n(n-1)}q^{\frac16n(n-1)(2n-1)+\frac{1}{2}n(n-1)t}\left\{\frac{(aq;q)_{t}}{(abq^2;q)_{t}}\right\}^{n}
  \nonumber\\&\times
  \prod_{k=1}^{n}\frac{(q,aq^{t+1},bq;q)_{n-k}}
  {(abq^{n-k+t+1};q)_{n-k}(abq^{t+2};q)_{2(n-k)}}.
  \label{eq:q-general-hankel}
  \end{align}
\end{corollary}
\begin{demo}{Proof}
If we use
\[
\mu_{n+t}
=\frac{(aq;q)_{n+t}}{(abq^2;q)_{n+t}}
=\frac{(aq;q)_{t}}{(abq^2;q)_{t}}\cdot\frac{(aq^{t+1};q)_{n}}{(abq^{t+2};q)_{n}},
\]
then we have
\begin{align*}
\det\left(\mu_{i+j+t}\right)_{0\leq i,j\leq n-1}
&=\det\left(
\frac{(aq;q)_{t}}{(abq^2;q)_{t}}\cdot\frac{(aq^{t+1};q)_{i+j}}{(abq^{t+2};q)_{i+j}}
\right)_{0\leq i,j\leq n-1}
\\&
=\left\{\frac{(aq;q)_{t}}{(abq^2;q)_{t}}\right\}^{n}
\det\left(
\frac{(a'q;q)_{i+j}}{(a'bq^{2};q)_{i+j}}
\right)_{0\leq i,j\leq n-1},
\end{align*}
where $a'=aq^{t}$.
If we use \eqref{eq:q-catalan-det},
then we obtain \eqref{eq:q-general-hankel} by a straightforward computation.
\end{demo}

We can prove \eqref{eq:Catalan-det}, \eqref{eq:B-Catalan} and \eqref{eq:D-Catalan}
as a corollary of Corollary~\ref{cor:q-Catalan-Hankel}.
\begin{demo}{Proof of \eqref{eq:Catalan-det}, \eqref{eq:B-Catalan} and \eqref{eq:D-Catalan}}
If we substitute $a=q^{\alpha}$ and $b=q^{\beta}$ into $\nu_{n}$,
and we put $q\rightarrow1$,
then we obtain 
$
\mu_{n}\rightarrow\frac{(\alpha+1)_{n}}{(\alpha+\beta+2)_{n}}
$,
which we write $\nu_{n}$.
Thus
\eqref{eq:q-general-hankel},
leads to 
\begin{align*}
  \det\left(\nu_{i+j+t}\right)_{0\leq i,j\leq n-1}
  =
\nu_{t}^{n}
  \prod_{k=1}^{n}\frac{(n-k)!(\alpha+t+1)_{n-k}(\beta+1)_{n-k}}
  {(\alpha+\beta+t+n-k+1)_{n-k}(\alpha+\beta+t+2)_{2(n-k)}}.
\end{align*}
Note that
\[
\nu_{n}=\begin{cases}
C_{n}/2^{2n}
&\text{ if $\alpha=-\frac12$ and $\beta=\frac12$,}\\
B_{n}/2^{2n}
&\text{ if $\alpha=\frac12$ and $\beta=-\frac12$,}\\
D_{n}/2^{2n}
&\text{ if $\alpha=-\frac12$ and $\beta=-\frac12$.}
\end{cases}
\]
Hence we obtain
\[
2^{2n(n+t-1)}\det\left(\nu_{i+j+t}\right)_{0\leq i,j\leq n-1}=\begin{cases}
\det\left(C_{i+j+t}\right)_{0\leq i,j\leq n-1}
&\text{ if $\alpha=-\frac12$ and $\beta=\frac12$,}\\
\det\left(B_{i+j+t}\right)_{0\leq i,j\leq n-1}
&\text{ if $\alpha=\frac12$ and $\beta=-\frac12$,}\\
\det\left(D_{i+j+t}\right)_{0\leq i,j\leq n-1}
&\text{ if $\alpha=-\frac12$ and $\beta=-\frac12$.}
\end{cases}
\]
Thus we can prove \eqref{eq:Catalan-det}, \eqref{eq:B-Catalan} and \eqref{eq:D-Catalan}
by direct computations from the above identity.
\end{demo}

In fact we can also obtain the following generalization of \eqref{eq:Krattenthaler}.
\begin{theorem}\label{th:q-kratt}
Let $n$ be a positive integer,
and $k_{0}$, $\dots$, $k_{n-1}$ nonnegative integers.
Then we have
\begin{equation}
  \det\left(\mu_{k_{i}+j}\right)_{0\leq i,j\leq n-1}
  =a^{\binom{n}2}q^{\binom{n+1}3}
  \prod_{i=0}^{n-1}\frac{(aq;q)_{k_{i}}}{(abq^{2};q)_{k_{i}+n-1}}
  \prod_{0\leq i<j\leq n-1}(q^{k_{i}}-q^{k_{j}})
  \prod_{i=0}^{n-1}(bq;q)_{i}.
\label{eq:q-kratt}
\end{equation}
\end{theorem}



\section{Non-intersecting lattice paths}\label{sec:lattice}

In this section we give our first proof of Theorem~\ref{th:q-Catalan-Hankel}
using non-intersecting lattice paths.

Let $m$ and $n$ be nonnegative integers.
A \newterm{Dyck path} is, by definition, a lattice path
in the plane lattice $\Bbb{Z}^2$ consisting of two types of steps:
rise vector $(1,1)$ and fall vector $(1,-1)$,
which never passes below the $x$-axis.
We say a rise vector (resp. fall vector) whose origin is $(x,y)$ and ends at $(x+1,y+1)$ (resp. $(x+1,y-1)$)
has \newterm{height} $y$.
%
%
%
%
%
\begin{figure}[ht]
\begin{center}
\setlength{\unitlength}{0.5mm}
\begin{picture}(150,80)
%
%
\put( 10, 10){\vector( 1, 0){120}}
\put( 10, 10){\vector( 0, 1){ 60}}
\put( 10, 10){\vector( 1, 1){10}}
\put( 20, 20){\vector( 1,-1){10}}
\put( 30, 10){\vector( 1, 1){10}}
\put( 40, 20){\vector( 1, 1){10}}
\put( 50, 30){\vector( 1, 1){10}}
\put( 60, 40){\vector( 1,-1){10}}
\put( 70, 30){\vector( 1, 1){10}}
\put( 80, 40){\vector( 1,-1){10}}
%
%
\put( 10, 10){\circle*{2}}
\put( 90, 30){\circle*{2}}
\put( 20, 10){\circle*{1.0}}
\put( 30, 10){\circle*{1.0}}
\put( 40, 10){\circle*{1.0}}
\put( 50, 10){\circle*{1.0}}
\put( 60, 10){\circle*{1.0}}
\put( 70, 10){\circle*{1.0}}
\put( 80, 10){\circle*{1.0}}
\put( 90, 10){\circle*{1.0}}
\put(100, 10){\circle*{1.0}}
\put(110, 10){\circle*{1.0}}
\put(120, 10){\circle*{1.0}}
\put( 10, 20){\circle*{1.0}}
\put( 20, 20){\circle*{1.0}}
\put( 30, 20){\circle*{1.0}}
\put( 40, 20){\circle*{1.0}}
\put( 50, 20){\circle*{1.0}}
\put( 60, 20){\circle*{1.0}}
\put( 70, 20){\circle*{1.0}}
\put( 80, 20){\circle*{1.0}}
\put( 90, 20){\circle*{1.0}}
\put(100, 20){\circle*{1.0}}
\put(110, 20){\circle*{1.0}}
\put(120, 20){\circle*{1.0}}
\put( 10, 30){\circle*{1.0}}
\put( 20, 30){\circle*{1.0}}
\put( 30, 30){\circle*{1.0}}
\put( 40, 30){\circle*{1.0}}
\put( 50, 30){\circle*{1.0}}
\put( 60, 30){\circle*{1.0}}
\put( 70, 30){\circle*{1.0}}
\put( 80, 30){\circle*{1.0}}
\put( 90, 30){\circle*{1.0}}
\put(100, 30){\circle*{1.0}}
\put(110, 30){\circle*{1.0}}
\put(120, 30){\circle*{1.0}}
\put( 10, 40){\circle*{1.0}}
\put( 20, 40){\circle*{1.0}}
\put( 30, 40){\circle*{1.0}}
\put( 40, 40){\circle*{1.0}}
\put( 50, 40){\circle*{1.0}}
\put( 60, 40){\circle*{1.0}}
\put( 70, 40){\circle*{1.0}}
\put( 80, 40){\circle*{1.0}}
\put( 90, 40){\circle*{1.0}}
\put(100, 40){\circle*{1.0}}
\put(110, 40){\circle*{1.0}}
\put(120, 40){\circle*{1.0}}
\put( 10, 50){\circle*{1.0}}
\put( 20, 50){\circle*{1.0}}
\put( 30, 50){\circle*{1.0}}
\put( 40, 50){\circle*{1.0}}
\put( 50, 50){\circle*{1.0}}
\put( 60, 50){\circle*{1.0}}
\put( 70, 50){\circle*{1.0}}
\put( 80, 50){\circle*{1.0}}
\put( 90, 50){\circle*{1.0}}
\put(100, 50){\circle*{1.0}}
\put(110, 50){\circle*{1.0}}
\put(120, 50){\circle*{1.0}}
\put( 10, 60){\circle*{1.0}}
\put( 20, 60){\circle*{1.0}}
\put( 30, 60){\circle*{1.0}}
\put( 40, 60){\circle*{1.0}}
\put( 50, 60){\circle*{1.0}}
\put( 60, 60){\circle*{1.0}}
\put( 70, 60){\circle*{1.0}}
\put( 80, 60){\circle*{1.0}}
\put( 90, 60){\circle*{1.0}}
\put(100, 60){\circle*{1.0}}
\put(110, 60){\circle*{1.0}}
\put(120, 60){\circle*{1.0}}
\put(  5,  5){\makebox(0,0){$\scriptstyle(0,0)$}}
\put( 95, 25){\makebox(0,0){$\scriptstyle(8,2)$}}
%
\put(130,  5){\makebox(0,0){$\scriptstyle x$}}
\put(  5, 70){\makebox(0,0){$\scriptstyle y$}}
\put( 15, 18){\makebox(0,0){\footnotesize\textcolor{red}{$\scriptstyle0$}}}
\put( 25, 18){\makebox(0,0){\footnotesize\textcolor{red}{$\scriptstyle1$}}}
\put( 35, 18){\makebox(0,0){\footnotesize\textcolor{red}{$\scriptstyle0$}}}
\put( 45, 28){\makebox(0,0){\footnotesize\textcolor{red}{$\scriptstyle1$}}}
\put( 55, 38){\makebox(0,0){\footnotesize\textcolor{red}{$\scriptstyle2$}}}
\put( 65, 38){\makebox(0,0){\footnotesize\textcolor{red}{$\scriptstyle3$}}}
\put( 75, 38){\makebox(0,0){\footnotesize\textcolor{red}{$\scriptstyle2$}}}
\put( 85, 38){\makebox(0,0){\footnotesize\textcolor{red}{$\scriptstyle3$}}}
%
\end{picture}
\caption{A Dyck Path starting from $(0,0)$ and ending at $(8,2)$}\label{fig:DPath}
\end{center}
\end{figure}
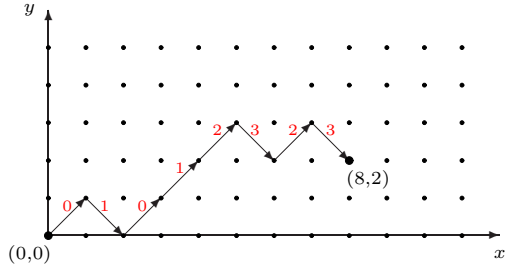
For example,
Figure~\ref{fig:DPath} presents a Dyck path starting from $(0,0)$ and ending at $(8,2)$,
in which each red number stands for the height of the step.
Let $\Cal{D}_{m,n}$ denote the set of Dyck paths starting from $(0,0)$ and ending at $(m,n)$.
Especially,
the cardinality of $\Cal{D}_{2n,0}$ is known to be the Catalan number $C_{n}$.

A \newterm{Motzkin path} is, by definition, a lattice path
in $\Bbb{Z}^2$ consisting of three types of steps:
rise vectors $(1,1)$, fall vectors $(1,-1)$, and (short) level vectors $(1,0)$
which never passes below the $x$-axis.
We say a rise vector, fall vector and level vector 
whose origin is $(x,y)$ and ends at $(x+1,y+1)$, $(x+1,y-1)$ and $(x+1,y)$
has \newterm{height} $y$, respectively.
%
%
%
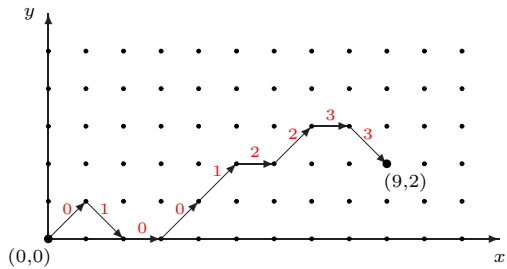
\begin{figure}[ht]
\begin{center}
\setlength{\unitlength}{0.5mm}
\begin{picture}(150,80)
%
%
\put( 10, 10){\vector( 1, 0){120}}
\put( 10, 10){\vector( 0, 1){ 60}}
\put( 10, 10){\vector( 1, 1){10}}
\put( 20, 20){\vector( 1,-1){10}}
\put( 30, 10){\vector( 1, 0){10}}
\put( 40, 10){\vector( 1, 1){10}}
\put( 50, 20){\vector( 1, 1){10}}
\put( 60, 30){\vector( 1, 0){10}}
\put( 70, 30){\vector( 1, 1){10}}
\put( 80, 40){\vector( 1, 0){10}}
\put( 90, 40){\vector( 1,-1){10}}
%
%
\put( 10, 10){\circle*{2}}
\put(100, 30){\circle*{2}}
\put( 20, 10){\circle*{1.0}}
\put( 30, 10){\circle*{1.0}}
\put( 40, 10){\circle*{1.0}}
\put( 50, 10){\circle*{1.0}}
\put( 60, 10){\circle*{1.0}}
\put( 70, 10){\circle*{1.0}}
\put( 80, 10){\circle*{1.0}}
\put( 90, 10){\circle*{1.0}}
\put(100, 10){\circle*{1.0}}
\put(110, 10){\circle*{1.0}}
\put(120, 10){\circle*{1.0}}
\put( 10, 20){\circle*{1.0}}
\put( 20, 20){\circle*{1.0}}
\put( 30, 20){\circle*{1.0}}
\put( 40, 20){\circle*{1.0}}
\put( 50, 20){\circle*{1.0}}
\put( 60, 20){\circle*{1.0}}
\put( 70, 20){\circle*{1.0}}
\put( 80, 20){\circle*{1.0}}
\put( 90, 20){\circle*{1.0}}
\put(100, 20){\circle*{1.0}}
\put(110, 20){\circle*{1.0}}
\put(120, 20){\circle*{1.0}}
\put( 10, 30){\circle*{1.0}}
\put( 20, 30){\circle*{1.0}}
\put( 30, 30){\circle*{1.0}}
\put( 40, 30){\circle*{1.0}}
\put( 50, 30){\circle*{1.0}}
\put( 60, 30){\circle*{1.0}}
\put( 70, 30){\circle*{1.0}}
\put( 80, 30){\circle*{1.0}}
\put( 90, 30){\circle*{1.0}}
\put(100, 30){\circle*{1.0}}
\put(110, 30){\circle*{1.0}}
\put(120, 30){\circle*{1.0}}
\put( 10, 40){\circle*{1.0}}
\put( 20, 40){\circle*{1.0}}
\put( 30, 40){\circle*{1.0}}
\put( 40, 40){\circle*{1.0}}
\put( 50, 40){\circle*{1.0}}
\put( 60, 40){\circle*{1.0}}
\put( 70, 40){\circle*{1.0}}
\put( 80, 40){\circle*{1.0}}
\put( 90, 40){\circle*{1.0}}
\put(100, 40){\circle*{1.0}}
\put(110, 40){\circle*{1.0}}
\put(120, 40){\circle*{1.0}}
\put( 10, 50){\circle*{1.0}}
\put( 20, 50){\circle*{1.0}}
\put( 30, 50){\circle*{1.0}}
\put( 40, 50){\circle*{1.0}}
\put( 50, 50){\circle*{1.0}}
\put( 60, 50){\circle*{1.0}}
\put( 70, 50){\circle*{1.0}}
\put( 80, 50){\circle*{1.0}}
\put( 90, 50){\circle*{1.0}}
\put(100, 50){\circle*{1.0}}
\put(110, 50){\circle*{1.0}}
\put(120, 50){\circle*{1.0}}
\put( 10, 60){\circle*{1.0}}
\put( 20, 60){\circle*{1.0}}
\put( 30, 60){\circle*{1.0}}
\put( 40, 60){\circle*{1.0}}
\put( 50, 60){\circle*{1.0}}
\put( 60, 60){\circle*{1.0}}
\put( 70, 60){\circle*{1.0}}
\put( 80, 60){\circle*{1.0}}
\put( 90, 60){\circle*{1.0}}
\put(100, 60){\circle*{1.0}}
\put(110, 60){\circle*{1.0}}
\put(120, 60){\circle*{1.0}}
\put(  5,  5){\makebox(0,0){$\scriptstyle(0,0)$}}
\put(105, 25){\makebox(0,0){$\scriptstyle(9,2)$}}
%
\put(130,  5){\makebox(0,0){$\scriptstyle x$}}
\put(  5, 70){\makebox(0,0){$\scriptstyle y$}}
\put( 15, 18){\makebox(0,0){\footnotesize\textcolor{red}{$\scriptstyle0$}}}
\put( 25, 18){\makebox(0,0){\footnotesize\textcolor{red}{$\scriptstyle1$}}}
\put( 35, 13){\makebox(0,0){\footnotesize\textcolor{red}{$\scriptstyle0$}}}
\put( 45, 18){\makebox(0,0){\footnotesize\textcolor{red}{$\scriptstyle0$}}}
\put( 55, 28){\makebox(0,0){\footnotesize\textcolor{red}{$\scriptstyle1$}}}
\put( 65, 33){\makebox(0,0){\footnotesize\textcolor{red}{$\scriptstyle2$}}}
\put( 75, 38){\makebox(0,0){\footnotesize\textcolor{red}{$\scriptstyle2$}}}
\put( 85, 43){\makebox(0,0){\footnotesize\textcolor{red}{$\scriptstyle3$}}}
\put( 95, 38){\makebox(0,0){\footnotesize\textcolor{red}{$\scriptstyle3$}}}
%
\end{picture}
\caption{A Moztkin path starting from $(0,0)$ and ending at $(9,2)$}\label{fig:MPath}
\end{center}
\end{figure}
Figure~\ref{fig:MPath} presents a Motzkin path starting from $(0,0)$ and ending at $(9,2)$,
in which each red number stands for the height of the step.
Let $\Cal{M}_{m,n}$ denote the set of Motzkin paths starting from $(0,0)$ and ending at $(m,n)$.
Note that 
the cardinality of $\Cal{M}_{n,0}$ is known to be the Motzkin number $M_{n}$.
We define the height of each step similarly as before.


%
%
%
%
%
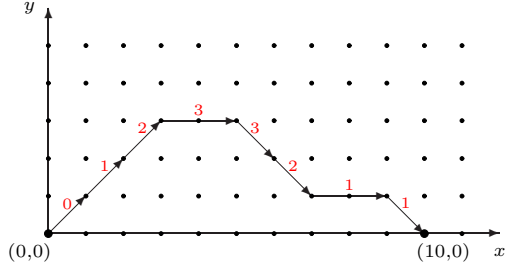
\begin{figure}[hbt]
\begin{center}
\setlength{\unitlength}{0.5mm}
\begin{picture}(150,80)
%
%
\put( 10, 10){\vector( 1, 0){120}}
\put( 10, 10){\vector( 0, 1){ 60}}
\put( 10, 10){\vector( 1, 1){10}}
\put( 20, 20){\vector( 1, 1){10}}
\put( 30, 30){\vector( 1, 1){10}}
\put( 40, 40){\vector( 1, 0){20}}
\put( 60, 40){\vector( 1,-1){10}}
\put( 70, 30){\vector( 1,-1){10}}
\put( 80, 20){\vector( 1, 0){20}}
\put(100, 20){\vector( 1,-1){10}}
%
%
\put( 10, 10){\circle*{2}}
\put(110, 10){\circle*{2}}
\put( 20, 10){\circle*{1.0}}
\put( 30, 10){\circle*{1.0}}
\put( 40, 10){\circle*{1.0}}
\put( 50, 10){\circle*{1.0}}
\put( 60, 10){\circle*{1.0}}
\put( 70, 10){\circle*{1.0}}
\put( 80, 10){\circle*{1.0}}
\put( 90, 10){\circle*{1.0}}
\put(100, 10){\circle*{1.0}}
\put(110, 10){\circle*{1.0}}
\put(120, 10){\circle*{1.0}}
\put( 10, 20){\circle*{1.0}}
\put( 20, 20){\circle*{1.0}}
\put( 30, 20){\circle*{1.0}}
\put( 40, 20){\circle*{1.0}}
\put( 50, 20){\circle*{1.0}}
\put( 60, 20){\circle*{1.0}}
\put( 70, 20){\circle*{1.0}}
\put( 80, 20){\circle*{1.0}}
\put( 90, 20){\circle*{1.0}}
\put(100, 20){\circle*{1.0}}
\put(110, 20){\circle*{1.0}}
\put(120, 20){\circle*{1.0}}
\put( 10, 30){\circle*{1.0}}
\put( 20, 30){\circle*{1.0}}
\put( 30, 30){\circle*{1.0}}
\put( 40, 30){\circle*{1.0}}
\put( 50, 30){\circle*{1.0}}
\put( 60, 30){\circle*{1.0}}
\put( 70, 30){\circle*{1.0}}
\put( 80, 30){\circle*{1.0}}
\put( 90, 30){\circle*{1.0}}
\put(100, 30){\circle*{1.0}}
\put(110, 30){\circle*{1.0}}
\put(120, 30){\circle*{1.0}}
\put( 10, 40){\circle*{1.0}}
\put( 20, 40){\circle*{1.0}}
\put( 30, 40){\circle*{1.0}}
\put( 40, 40){\circle*{1.0}}
\put( 50, 40){\circle*{1.0}}
\put( 60, 40){\circle*{1.0}}
\put( 70, 40){\circle*{1.0}}
\put( 80, 40){\circle*{1.0}}
\put( 90, 40){\circle*{1.0}}
\put(100, 40){\circle*{1.0}}
\put(110, 40){\circle*{1.0}}
\put(120, 40){\circle*{1.0}}
\put( 10, 50){\circle*{1.0}}
\put( 20, 50){\circle*{1.0}}
\put( 30, 50){\circle*{1.0}}
\put( 40, 50){\circle*{1.0}}
\put( 50, 50){\circle*{1.0}}
\put( 60, 50){\circle*{1.0}}
\put( 70, 50){\circle*{1.0}}
\put( 80, 50){\circle*{1.0}}
\put( 90, 50){\circle*{1.0}}
\put(100, 50){\circle*{1.0}}
\put(110, 50){\circle*{1.0}}
\put(120, 50){\circle*{1.0}}
\put( 10, 60){\circle*{1.0}}
\put( 20, 60){\circle*{1.0}}
\put( 30, 60){\circle*{1.0}}
\put( 40, 60){\circle*{1.0}}
\put( 50, 60){\circle*{1.0}}
\put( 60, 60){\circle*{1.0}}
\put( 70, 60){\circle*{1.0}}
\put( 80, 60){\circle*{1.0}}
\put( 90, 60){\circle*{1.0}}
\put(100, 60){\circle*{1.0}}
\put(110, 60){\circle*{1.0}}
\put(120, 60){\circle*{1.0}}
\put(  5,  5){\makebox(0,0){$\scriptstyle(0,0)$}}
\put(115,  5){\makebox(0,0){$\scriptstyle(10,0)$}}
%
\put(130,  5){\makebox(0,0){$\scriptstyle x$}}
\put(  5, 70){\makebox(0,0){$\scriptstyle y$}}
\put( 15, 18){\makebox(0,0){\footnotesize\textcolor{red}{$\scriptstyle0$}}}
\put( 25, 28){\makebox(0,0){\footnotesize\textcolor{red}{$\scriptstyle1$}}}
\put( 35, 38){\makebox(0,0){\footnotesize\textcolor{red}{$\scriptstyle2$}}}
\put( 50, 43){\makebox(0,0){\footnotesize\textcolor{red}{$\scriptstyle3$}}}
\put( 65, 38){\makebox(0,0){\footnotesize\textcolor{red}{$\scriptstyle3$}}}
\put( 75, 28){\makebox(0,0){\footnotesize\textcolor{red}{$\scriptstyle2$}}}
\put( 90, 23){\makebox(0,0){\footnotesize\textcolor{red}{$\scriptstyle1$}}}
\put(105, 18){\makebox(0,0){\footnotesize\textcolor{red}{$\scriptstyle1$}}}
%
\end{picture}
\caption{A Schr\"oder path starting from $(0,0)$ and ending at $(10,0)$}\label{fig:SPath}
\end{center}
\end{figure}
A \newterm{Schr\"oder path} is, by definition, a lattice path
in $\Bbb{Z}^2$ consisting of three types of steps:
rise vectors $(1,1)$, fall vectors $(1,-1)$, and long level vectors $(2,0)$
which never passes below the $x$-axis.
Figure~\ref{fig:SPath} presents a Schr\"oder path starting from $(0,0)$ and ending at $(10,0)$,
in which each red number stands for the height of the step.
Let $\Cal{S}_{m,n}$ denote the set of Schr\"oder paths starting from $(0,0)$ and ending at $(m,n)$.
Note that
the cardinality of $\Cal{S}_{2n,0}$ is known to be the large Schr\"oder number $S_{n}$.

Assign the weight $a_{h}$, $b_{h}$, $c_{h}$ to each rise vector,
fall vector, (short or long) level vector of height $h$, 
respectively.
Set the weight of a path $P$ to be the product of the weights of its edges and denote it by $w(P)$.
Given any family $\mathcal{F}$ of paths,
we write the generating function of $\mathcal{F}$ as
\[
\GF{\mathcal{F}}=\sum_{P\in\mathcal{F}}w(P).
\]
\begin{prop}(Flajolet \cite{F})
The generating function for the Dyck paths
is given by the following
Stieltjes type continued fraction:
\begin{align*}
&\sum_{n\geq0}\GF{\Cal{D}_{(2n,0)}}t^{2n}
=\frac1{1-\frac{a_0b_1t^2}{1-\frac{a_1b_2t^2}{1-\frac{a_2b_3t^2}{\ddots}}}}.
\end{align*}
Meanwhile, the generating function for the Motzkin paths
is given by the following
Jacobi type continued fraction:
\begin{align*}
&\sum_{n\geq0}\GF{\Cal{M}_{(n,0)}}t^{n}
=\frac1{1-c_0t-\frac{a_0b_1t^2}{1-c_1t-\frac{a_1b_2t^2}{1-c_2t-\frac{a_2b_3t^2}{\ddots}}}}.
\end{align*}
\end{prop}
It is also easy to see the following proposition holds.
\begin{prop}
Let $n$ be a positive integer.
Then the generating function for 
Schr\"oder paths
is given by the following
continued fraction:
\begin{align*}
\sum_{n\geq0}\GF{\Cal{S}_{(2n,0)}}t^{2n}
=\frac1{1-c_0t^2-\frac{a_0b_1t^2}{1-c_1t^2-\frac{a_1b_2t^2}{1-c_2t^2-\frac{a_2b_3t^2}{\ddots}}}}.
\end{align*}
\end{prop}

\bigskip

Next we recall notation and definitions used for the lattice path method due to Gessel and Viennot \cite{GV}.
Let $D=(V,E)$ be an acyclic digraph without multiple edges.
If $u$ and $v$ are any pair of vertices,
let $\PATH{u}{v}$ denote the set of all directed paths from $u$ to $v$.
For a fixed positive integer $n$,
an \defterm{$n$-vertex} is an $n$-tuple of vertices of $D$.
If $\boldsymbol{u}=(u_1,\dots,u_n)$ and $\boldsymbol{v}=(v_1,\dots,v_n)$ are $n$-vertices,
an \defterm{$n$-path} from $\boldsymbol{u}$ to $\boldsymbol{v}$ is an $n$-tuple $\boldsymbol{P}=(P_1,\dots,P_n)$
such that $P_i\in\PATH{u_i}{v_i}$, $i=1,\dots,n$.
The $n$-path $\boldsymbol{P}=(P_1,\dots,P_n)$ is said to be \defterm{non-intersecting}
if any two different paths $P_i$ and $P_j$ have no vertex in common.
We will write $\PATH{\boldsymbol{u}}{\boldsymbol{v}}$ for the set of all $n$-paths from $\boldsymbol{u}$ to $\boldsymbol{v}$,
and write $\NPATH{\boldsymbol{u}}{\boldsymbol{v}}$ for the subset of $\PATH{\boldsymbol{u}}{\boldsymbol{v}}$ 
consisting of non-intersecting $n$-paths.
If $\boldsymbol{u}=(u_1,\dots,u_m)$ and $\boldsymbol{v}=(v_1,\dots,v_n)$ are linearly ordered sets of vertices of $D$,
then $\boldsymbol{u}$ is said to be \defterm{$D$-compatible} with $\boldsymbol{v}$ if
every path $P\in{\cal P}(u_i,v_l)$ intersects with every path $Q\in{\cal P}(u_j,v_k)$ whenever $i<j$ and $k<l$.
Let $S_n$ denote the symmetric group on $\{1,2,\dots,n\}$.
Then for $\pi\in S_n$,
by $\boldsymbol{v}^\pi$ we mean the $n$ vertex $(v_{\pi(1)},\dots,v_{\pi(n)})$.

The weight $w(\boldsymbol{P})$ of an $n$-path $\boldsymbol{P}$ is defined to be the product
of the weights of its components.
Thus,
if $\boldsymbol{u}=(u_1,\dots,u_n)$ and $\boldsymbol{v}=(v_1,\dots,v_n)$ are $n$-vertices,
 we define the generating functions $F(\boldsymbol{u},\boldsymbol{v})=\GF{\PATH{\boldsymbol{u}}{\boldsymbol{v}}}=\sum_{\boldsymbol{P}\in\PATH{\boldsymbol{u}}{\boldsymbol{v}}}w(\boldsymbol{P})$ 
and $F_0(\boldsymbol{u},\boldsymbol{v})=\GF{\NPATH{\boldsymbol{u}}{\boldsymbol{v}}}=\sum_{\boldsymbol{P}\in\NPATH{\boldsymbol{u}}{\boldsymbol{v}}}w(\boldsymbol{P})$.
In particular, if $u$ and $v$ are any pair of vertices,
we write
\begin{equation*}
h(u,v)=\GF{\PATH{u}{v}}=\sum_{P\in\PATH{u}{v}}w(P).
\end{equation*}
The following lemma is called the Gessel-Viennot formula for counting lattice paths
in terms of determinants.
(See \cite{GV}.)
\begin{lemma}
\label{lem:Gessel-Viennot}
(Lidstr\"om-Gessel-Viennot)

Let $\boldsymbol{u}=(u_1,\dots,u_n)$ and $\boldsymbol{v}=(v_1,\dots,v_n)$ be two $n$-vertices in an acyclic digraph $D$.
Then
\begin{equation}
\sum_{\pi\in S_n}\sgn\pi\ F_0(\boldsymbol{u}^{\pi},\boldsymbol{v})=\det[h(u_i,v_j)]_{1\le i,j\le n}.
\label{eq:LGV1}
\end{equation}
In particular, if $\boldsymbol{u}$ is $D$-compatible with $\boldsymbol{v}$,
then
\begin{equation}
F_0(\boldsymbol{u},\boldsymbol{v})=\det[h(u_i,v_j)]_{1\le i,j\le n}.
\label{eq:LGV2}
\end{equation}
\end{lemma}
If we apply Lemma~\ref{lem:Gessel-Viennot} to Dyck paths,
then we obtain the following proposition:
\begin{prop}
\label{prop:viennot}
Let $G_{m}=\GF{\Cal{D}_{2m,0}}$ for non-negative integer $m$.
\par\noindent
(i)
If $t=0$,
then we have
\begin{equation}
\det\left(G_{i+j}\right)_{0\leq i,j\leq n-1}
=\prod_{i=1}^{n}\left(a_{2i-2}b_{2i-1}a_{2i-1}b_{2i}\right)^{n-i}.
\label{eq:t=0}
\end{equation}
(ii)
If $t=1$,
then we have
\begin{equation}
\det\left(G_{i+j+1}\right)_{0\leq i,j\leq n-1}
=\prod_{i=1}^{n}\left(a_{2i-2}b_{2i-1}\right)^{n-i+1}\left(a_{2i-1}b_{2i}\right)^{n-i}.
\label{eq:t=1}
\end{equation}
(iii)
If $t=2$,
then we have
$\det\left(G_{i+j+2}\right)_{0\leq i,j\leq n-1}$ equals
\begin{align}
\sum_{k=0}^{n}\prod_{i=1}^{k}
\left(a_{0}a_{1}\cdots a_{2i-3}a_{2i-2}^2b_{1}b_{2}\cdots b_{2i-1}b_{2i-1}^2\right)
\cdot\prod_{i=1}^{k}\left(a_{0}a_{1}\cdots a_{2i-1}b_{1}b_{2}\cdots b_{2i}\right).
\label{eq:t=2}
\end{align}
(iv)
If $t=3$,
then we have $\det\left(G_{i+j+3}\right)_{0\leq i,j\leq n-1}$ equals
\begin{align}
&\sum_{k=0}^{n}
\Bigl\{\sum_{l=0}^{k}
\prod_{i=1}^{l}\left(a_{0}a_{1}\cdots a_{2i-3}a_{2i-2}^2b_{2i-1}\right)
\prod_{i=l+1}^{k}\left(a_{0}a_{1}\cdots a_{2i-3}a_{2i-2}a_{2i-1}b_{2i}\right)\Bigr\}
\nonumber\\
&\times\Bigl\{\sum_{l=0}^{k}
\prod_{i=1}^{l}\left(b_{1}b_{2}\cdots b_{2i-2}b_{2i-1}^2a_{2i-2}\right)
\prod_{i=l+1}^{k}\left(b_{1}b_{2}\cdots b_{2i-2}b_{2i-1}b_{2i}a_{2i-1}\right)\Bigr\}
\nonumber\\&\times
\prod_{i=k+1}^{n}\left(a_{0}a_{1}\cdots a_{2i-1}a_{2i}b_{1}b_{2}\cdots b_{2i}b_{2i+1}\right).
\label{eq:t=3}
\end{align}
((i) and (ii) of this proposition are originally appeared in \cite[Ch.~4, \S3]{V0}.)
\end{prop}
\begin{demo}{Proof}
We consider the digraph $(V,E)$,
in which $V$ is the plane lattice $\Bbb{Z}^2$
and $E$ the set of rise vectors and fall vectors in the above half plane.
Let $u_{i}=(x_{0}-2(i-1),0)$ and $v_{j}=(x_{0}+2(j+t-1),0)$
for $i,j=1,2,\dots,n$, $t=0,1,2,3$ and a fixed integer $x_{0}$.
It is easy to see that the $n$-vertex 
$\boldsymbol{u}=(u_1,\dots,u_n)$ is \defterm{$D$-compatible} 
with the $n$-vertex $\boldsymbol{v}=(v_1,\dots,v_n)$.
%
%
%
\begin{figure}[ht]
\begin{center}
\setlength{\unitlength}{0.4mm}
\begin{picture}(150,90)
%
%
\put( 10, 10){\vector( 1, 0){140}}
\put( 10, 10){\vector( 0, 1){ 80}}
\put( 60, 10){\vector( 1, 1){10}}
\put( 70, 20){\vector( 1, 1){10}}
\put( 80, 30){\vector( 1,-1){10}}
\put( 90, 20){\vector( 1,-1){10}}
\put( 40, 10){\vector( 1, 1){10}}
\put( 50, 20){\vector( 1, 1){10}}
\put( 60, 30){\vector( 1, 1){10}}
\put( 70, 40){\vector( 1, 1){10}}
\put( 80, 50){\vector( 1,-1){10}}
\put( 90, 40){\vector( 1,-1){10}}
\put(100, 30){\vector( 1,-1){10}}
\put(110, 20){\vector( 1,-1){10}}
\put( 20, 10){\vector( 1, 1){10}}
\put( 30, 20){\vector( 1, 1){10}}
\put( 40, 30){\vector( 1, 1){10}}
\put( 50, 40){\vector( 1, 1){10}}
\put( 60, 50){\vector( 1, 1){10}}
\put( 70, 60){\vector( 1, 1){10}}
\put( 80, 70){\vector( 1,-1){10}}
\put( 90, 60){\vector( 1,-1){10}}
\put(100, 50){\vector( 1,-1){10}}
\put(110, 40){\vector( 1,-1){10}}
\put(120, 30){\vector( 1,-1){10}}
\put(130, 20){\vector( 1,-1){10}}
%
%
\put( 20, 10){\circle*{2}}
\put( 40, 10){\circle*{2}}
\put( 60, 10){\circle*{2}}
\put( 80, 10){\circle*{2}}
\put(100, 10){\circle*{2}}
\put(120, 10){\circle*{2}}
\put(140, 10){\circle*{2}}
\put( 10, 10){\circle*{1.0}}
\put( 20, 10){\circle*{1.0}}
\put( 30, 10){\circle*{1.0}}
\put( 40, 10){\circle*{1.0}}
\put( 50, 10){\circle*{1.0}}
\put( 60, 10){\circle*{1.0}}
\put( 70, 10){\circle*{1.0}}
\put( 80, 10){\circle*{1.0}}
\put( 90, 10){\circle*{1.0}}
\put(100, 10){\circle*{1.0}}
\put(110, 10){\circle*{1.0}}
\put(120, 10){\circle*{1.0}}
\put(130, 10){\circle*{1.0}}
\put(140, 10){\circle*{1.0}}
\put( 10, 20){\circle*{1.0}}
\put( 20, 20){\circle*{1.0}}
\put( 30, 20){\circle*{1.0}}
\put( 40, 20){\circle*{1.0}}
\put( 50, 20){\circle*{1.0}}
\put( 60, 20){\circle*{1.0}}
\put( 70, 20){\circle*{1.0}}
\put( 80, 20){\circle*{1.0}}
\put( 90, 20){\circle*{1.0}}
\put(100, 20){\circle*{1.0}}
\put(110, 20){\circle*{1.0}}
\put(120, 20){\circle*{1.0}}
\put(130, 20){\circle*{1.0}}
\put(140, 20){\circle*{1.0}}
\put( 10, 30){\circle*{1.0}}
\put( 20, 30){\circle*{1.0}}
\put( 30, 30){\circle*{1.0}}
\put( 40, 30){\circle*{1.0}}
\put( 50, 30){\circle*{1.0}}
\put( 60, 30){\circle*{1.0}}
\put( 70, 30){\circle*{1.0}}
\put( 80, 30){\circle*{1.0}}
\put( 90, 30){\circle*{1.0}}
\put(100, 30){\circle*{1.0}}
\put(110, 30){\circle*{1.0}}
\put(120, 30){\circle*{1.0}}
\put(130, 30){\circle*{1.0}}
\put(140, 30){\circle*{1.0}}
\put( 10, 40){\circle*{1.0}}
\put( 20, 40){\circle*{1.0}}
\put( 30, 40){\circle*{1.0}}
\put( 40, 40){\circle*{1.0}}
\put( 50, 40){\circle*{1.0}}
\put( 60, 40){\circle*{1.0}}
\put( 70, 40){\circle*{1.0}}
\put( 80, 40){\circle*{1.0}}
\put( 90, 40){\circle*{1.0}}
\put(100, 40){\circle*{1.0}}
\put(110, 40){\circle*{1.0}}
\put(120, 40){\circle*{1.0}}
\put(130, 40){\circle*{1.0}}
\put(140, 40){\circle*{1.0}}
\put( 10, 50){\circle*{1.0}}
\put( 20, 50){\circle*{1.0}}
\put( 30, 50){\circle*{1.0}}
\put( 40, 50){\circle*{1.0}}
\put( 50, 50){\circle*{1.0}}
\put( 60, 50){\circle*{1.0}}
\put( 70, 50){\circle*{1.0}}
\put( 80, 50){\circle*{1.0}}
\put( 90, 50){\circle*{1.0}}
\put(100, 50){\circle*{1.0}}
\put(110, 50){\circle*{1.0}}
\put(120, 50){\circle*{1.0}}
\put(130, 50){\circle*{1.0}}
\put(140, 50){\circle*{1.0}}
\put( 10, 60){\circle*{1.0}}
\put( 20, 60){\circle*{1.0}}
\put( 30, 60){\circle*{1.0}}
\put( 40, 60){\circle*{1.0}}
\put( 50, 60){\circle*{1.0}}
\put( 60, 60){\circle*{1.0}}
\put( 70, 60){\circle*{1.0}}
\put( 80, 60){\circle*{1.0}}
\put( 90, 60){\circle*{1.0}}
\put(100, 60){\circle*{1.0}}
\put(110, 60){\circle*{1.0}}
\put(120, 60){\circle*{1.0}}
\put(130, 60){\circle*{1.0}}
\put(140, 60){\circle*{1.0}}
\put( 10, 70){\circle*{1.0}}
\put( 20, 70){\circle*{1.0}}
\put( 30, 70){\circle*{1.0}}
\put( 40, 70){\circle*{1.0}}
\put( 50, 70){\circle*{1.0}}
\put( 60, 70){\circle*{1.0}}
\put( 70, 70){\circle*{1.0}}
\put( 80, 70){\circle*{1.0}}
\put( 90, 70){\circle*{1.0}}
\put(100, 70){\circle*{1.0}}
\put(110, 70){\circle*{1.0}}
\put(120, 70){\circle*{1.0}}
\put(130, 70){\circle*{1.0}}
\put(140, 70){\circle*{1.0}}
\put( 10, 80){\circle*{1.0}}
\put( 20, 80){\circle*{1.0}}
\put( 30, 80){\circle*{1.0}}
\put( 40, 80){\circle*{1.0}}
\put( 50, 80){\circle*{1.0}}
\put( 60, 80){\circle*{1.0}}
\put( 70, 80){\circle*{1.0}}
\put( 80, 80){\circle*{1.0}}
\put( 90, 80){\circle*{1.0}}
\put(100, 80){\circle*{1.0}}
\put(110, 80){\circle*{1.0}}
\put(120, 80){\circle*{1.0}}
\put(130, 80){\circle*{1.0}}
\put(140, 80){\circle*{1.0}}
\put(  5,  5){\makebox(0,0){$\scriptstyle(0,0)$}}
%
\put(150,  5){\makebox(0,0){$\scriptstyle x$}}
\put(  5, 90){\makebox(0,0){$\scriptstyle y$}}
\put( 20,  5){\makebox(0,0){\footnotesize\textcolor{blue}{$\scriptstyle u_4$}}}
\put( 40,  5){\makebox(0,0){\footnotesize\textcolor{blue}{$\scriptstyle u_3$}}}
\put( 60,  5){\makebox(0,0){\footnotesize\textcolor{blue}{$\scriptstyle u_2$}}}
\put( 80,  5){\makebox(0,0){\footnotesize\textcolor{blue}{$\scriptstyle u_1=v_1$}}}
\put(100,  5){\makebox(0,0){\footnotesize\textcolor{blue}{$\scriptstyle v_2$}}}
\put(120,  5){\makebox(0,0){\footnotesize\textcolor{blue}{$\scriptstyle v_3$}}}
\put(140,  5){\makebox(0,0){\footnotesize\textcolor{blue}{$\scriptstyle v_4$}}}
\put( 25, 18){\makebox(0,0){\footnotesize\textcolor{red}{$\scriptstyle a_0$}}}
\put( 45, 18){\makebox(0,0){\footnotesize\textcolor{red}{$\scriptstyle a_0$}}}
\put( 65, 18){\makebox(0,0){\footnotesize\textcolor{red}{$\scriptstyle a_0$}}}
\put( 95, 18){\makebox(0,0){\footnotesize\textcolor{red}{$\scriptstyle b_1$}}}
\put(115, 18){\makebox(0,0){\footnotesize\textcolor{red}{$\scriptstyle b_1$}}}
\put(135, 18){\makebox(0,0){\footnotesize\textcolor{red}{$\scriptstyle b_1$}}}
\put( 35, 28){\makebox(0,0){\footnotesize\textcolor{red}{$\scriptstyle a_1$}}}
\put( 55, 28){\makebox(0,0){\footnotesize\textcolor{red}{$\scriptstyle a_1$}}}
\put( 75, 28){\makebox(0,0){\footnotesize\textcolor{red}{$\scriptstyle a_1$}}}
\put( 85, 28){\makebox(0,0){\footnotesize\textcolor{red}{$\scriptstyle b_2$}}}
\put(105, 28){\makebox(0,0){\footnotesize\textcolor{red}{$\scriptstyle b_2$}}}
\put(125, 28){\makebox(0,0){\footnotesize\textcolor{red}{$\scriptstyle b_2$}}}
\put( 45, 38){\makebox(0,0){\footnotesize\textcolor{red}{$\scriptstyle a_2$}}}
\put( 65, 38){\makebox(0,0){\footnotesize\textcolor{red}{$\scriptstyle a_2$}}}
\put( 95, 38){\makebox(0,0){\footnotesize\textcolor{red}{$\scriptstyle b_3$}}}
\put(115, 38){\makebox(0,0){\footnotesize\textcolor{red}{$\scriptstyle b_3$}}}
\put( 55, 48){\makebox(0,0){\footnotesize\textcolor{red}{$\scriptstyle a_3$}}}
\put( 75, 48){\makebox(0,0){\footnotesize\textcolor{red}{$\scriptstyle a_3$}}}
\put( 85, 48){\makebox(0,0){\footnotesize\textcolor{red}{$\scriptstyle b_4$}}}
\put(105, 48){\makebox(0,0){\footnotesize\textcolor{red}{$\scriptstyle b_4$}}}
\put( 65, 58){\makebox(0,0){\footnotesize\textcolor{red}{$\scriptstyle a_4$}}}
\put( 95, 58){\makebox(0,0){\footnotesize\textcolor{red}{$\scriptstyle b_5$}}}
\put( 75, 68){\makebox(0,0){\footnotesize\textcolor{red}{$\scriptstyle a_5$}}}
\put( 85, 68){\makebox(0,0){\footnotesize\textcolor{red}{$\scriptstyle b_6$}}}
\end{picture}
\caption{$t=0$ and $n=4$}\label{fig:t=0}
\end{center}
\end{figure}
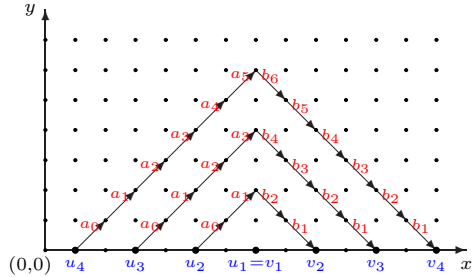
If $t=0$,
then there is always a unique $n$-path $\boldsymbol{P}=(P_1,\dots,P_n)$
that connect $\boldsymbol{u}$ to $\boldsymbol{v}$
as in Figure~\ref{fig:t=0}.
By multiplying the weights of all edges in $\boldsymbol{P}$,
we obtain the right-hand side of \eqref{eq:t=0}.
On the other hand,
applying Lemma~\ref{lem:Gessel-Viennot},
we  obtain the left-hand side of \eqref{eq:t=0}.

%
%
%
\begin{figure}[ht]
\begin{center}
\setlength{\unitlength}{0.4mm}
\begin{picture}(170,110)
%
%
\put( 10, 10){\vector( 1, 0){160}}
\put( 10, 10){\vector( 0, 1){100}}
\put( 80, 10){\vector( 1, 1){10}}
\put( 90, 20){\vector( 1,-1){10}}
\put( 60, 10){\vector( 1, 1){10}}
\put( 70, 20){\vector( 1, 1){10}}
\put( 80, 30){\vector( 1, 1){10}}
\put( 90, 40){\vector( 1,-1){10}}
\put(100, 30){\vector( 1,-1){10}}
\put(110, 20){\vector( 1,-1){10}}
\put( 40, 10){\vector( 1, 1){10}}
\put( 50, 20){\vector( 1, 1){10}}
\put( 60, 30){\vector( 1, 1){10}}
\put( 70, 40){\vector( 1, 1){10}}
\put( 80, 50){\vector( 1, 1){10}}
\put( 90, 60){\vector( 1,-1){10}}
\put(100, 50){\vector( 1,-1){10}}
\put(110, 40){\vector( 1,-1){10}}
\put(120, 30){\vector( 1,-1){10}}
\put(130, 20){\vector( 1,-1){10}}
\put( 20, 10){\vector( 1, 1){10}}
\put( 30, 20){\vector( 1, 1){10}}
\put( 40, 30){\vector( 1, 1){10}}
\put( 50, 40){\vector( 1, 1){10}}
\put( 60, 50){\vector( 1, 1){10}}
\put( 70, 60){\vector( 1, 1){10}}
\put( 80, 70){\vector( 1, 1){10}}
\put( 90, 80){\vector( 1,-1){10}}
\put(100, 70){\vector( 1,-1){10}}
\put(110, 60){\vector( 1,-1){10}}
\put(120, 50){\vector( 1,-1){10}}
\put(130, 40){\vector( 1,-1){10}}
\put(140, 30){\vector( 1,-1){10}}
\put(150, 20){\vector( 1,-1){10}}
%
%
\put( 20, 10){\circle*{2}}
\put( 40, 10){\circle*{2}}
\put( 60, 10){\circle*{2}}
\put( 80, 10){\circle*{2}}
\put(100, 10){\circle*{2}}
\put(120, 10){\circle*{2}}
\put(140, 10){\circle*{2}}
\put(160, 10){\circle*{2}}
\put( 10, 10){\circle*{1.0}}
\put( 20, 10){\circle*{1.0}}
\put( 30, 10){\circle*{1.0}}
\put( 40, 10){\circle*{1.0}}
\put( 50, 10){\circle*{1.0}}
\put( 60, 10){\circle*{1.0}}
\put( 70, 10){\circle*{1.0}}
\put( 80, 10){\circle*{1.0}}
\put( 90, 10){\circle*{1.0}}
\put(100, 10){\circle*{1.0}}
\put(110, 10){\circle*{1.0}}
\put(120, 10){\circle*{1.0}}
\put(130, 10){\circle*{1.0}}
\put(140, 10){\circle*{1.0}}
\put(150, 10){\circle*{1.0}}
\put(160, 10){\circle*{1.0}}
\put( 10, 20){\circle*{1.0}}
\put( 20, 20){\circle*{1.0}}
\put( 30, 20){\circle*{1.0}}
\put( 40, 20){\circle*{1.0}}
\put( 50, 20){\circle*{1.0}}
\put( 60, 20){\circle*{1.0}}
\put( 70, 20){\circle*{1.0}}
\put( 80, 20){\circle*{1.0}}
\put( 90, 20){\circle*{1.0}}
\put(100, 20){\circle*{1.0}}
\put(110, 20){\circle*{1.0}}
\put(120, 20){\circle*{1.0}}
\put(130, 20){\circle*{1.0}}
\put(140, 20){\circle*{1.0}}
\put(150, 20){\circle*{1.0}}
\put(160, 20){\circle*{1.0}}
\put( 10, 30){\circle*{1.0}}
\put( 20, 30){\circle*{1.0}}
\put( 30, 30){\circle*{1.0}}
\put( 40, 30){\circle*{1.0}}
\put( 50, 30){\circle*{1.0}}
\put( 60, 30){\circle*{1.0}}
\put( 70, 30){\circle*{1.0}}
\put( 80, 30){\circle*{1.0}}
\put( 90, 30){\circle*{1.0}}
\put(100, 30){\circle*{1.0}}
\put(110, 30){\circle*{1.0}}
\put(120, 30){\circle*{1.0}}
\put(130, 30){\circle*{1.0}}
\put(140, 30){\circle*{1.0}}
\put(150, 30){\circle*{1.0}}
\put(160, 30){\circle*{1.0}}
\put( 10, 40){\circle*{1.0}}
\put( 20, 40){\circle*{1.0}}
\put( 30, 40){\circle*{1.0}}
\put( 40, 40){\circle*{1.0}}
\put( 50, 40){\circle*{1.0}}
\put( 60, 40){\circle*{1.0}}
\put( 70, 40){\circle*{1.0}}
\put( 80, 40){\circle*{1.0}}
\put( 90, 40){\circle*{1.0}}
\put(100, 40){\circle*{1.0}}
\put(110, 40){\circle*{1.0}}
\put(120, 40){\circle*{1.0}}
\put(130, 40){\circle*{1.0}}
\put(140, 40){\circle*{1.0}}
\put(150, 40){\circle*{1.0}}
\put(160, 40){\circle*{1.0}}
\put( 10, 50){\circle*{1.0}}
\put( 20, 50){\circle*{1.0}}
\put( 30, 50){\circle*{1.0}}
\put( 40, 50){\circle*{1.0}}
\put( 50, 50){\circle*{1.0}}
\put( 60, 50){\circle*{1.0}}
\put( 70, 50){\circle*{1.0}}
\put( 80, 50){\circle*{1.0}}
\put( 90, 50){\circle*{1.0}}
\put(100, 50){\circle*{1.0}}
\put(110, 50){\circle*{1.0}}
\put(120, 50){\circle*{1.0}}
\put(130, 50){\circle*{1.0}}
\put(140, 50){\circle*{1.0}}
\put(150, 50){\circle*{1.0}}
\put(160, 50){\circle*{1.0}}
\put( 10, 60){\circle*{1.0}}
\put( 20, 60){\circle*{1.0}}
\put( 30, 60){\circle*{1.0}}
\put( 40, 60){\circle*{1.0}}
\put( 50, 60){\circle*{1.0}}
\put( 60, 60){\circle*{1.0}}
\put( 70, 60){\circle*{1.0}}
\put( 80, 60){\circle*{1.0}}
\put( 90, 60){\circle*{1.0}}
\put(100, 60){\circle*{1.0}}
\put(110, 60){\circle*{1.0}}
\put(120, 60){\circle*{1.0}}
\put(130, 60){\circle*{1.0}}
\put(140, 60){\circle*{1.0}}
\put(150, 60){\circle*{1.0}}
\put(160, 60){\circle*{1.0}}
\put( 10, 70){\circle*{1.0}}
\put( 20, 70){\circle*{1.0}}
\put( 30, 70){\circle*{1.0}}
\put( 40, 70){\circle*{1.0}}
\put( 50, 70){\circle*{1.0}}
\put( 60, 70){\circle*{1.0}}
\put( 70, 70){\circle*{1.0}}
\put( 80, 70){\circle*{1.0}}
\put( 90, 70){\circle*{1.0}}
\put(100, 70){\circle*{1.0}}
\put(110, 70){\circle*{1.0}}
\put(120, 70){\circle*{1.0}}
\put(130, 70){\circle*{1.0}}
\put(140, 70){\circle*{1.0}}
\put(150, 70){\circle*{1.0}}
\put(160, 70){\circle*{1.0}}
\put( 10, 80){\circle*{1.0}}
\put( 20, 80){\circle*{1.0}}
\put( 30, 80){\circle*{1.0}}
\put( 40, 80){\circle*{1.0}}
\put( 50, 80){\circle*{1.0}}
\put( 60, 80){\circle*{1.0}}
\put( 70, 80){\circle*{1.0}}
\put( 80, 80){\circle*{1.0}}
\put( 90, 80){\circle*{1.0}}
\put(100, 80){\circle*{1.0}}
\put(110, 80){\circle*{1.0}}
\put(120, 80){\circle*{1.0}}
\put(130, 80){\circle*{1.0}}
\put(140, 80){\circle*{1.0}}
\put(150, 80){\circle*{1.0}}
\put(160, 80){\circle*{1.0}}
\put( 10, 90){\circle*{1.0}}
\put( 20, 90){\circle*{1.0}}
\put( 30, 90){\circle*{1.0}}
\put( 40, 90){\circle*{1.0}}
\put( 50, 90){\circle*{1.0}}
\put( 60, 90){\circle*{1.0}}
\put( 70, 90){\circle*{1.0}}
\put( 80, 90){\circle*{1.0}}
\put( 90, 90){\circle*{1.0}}
\put(100, 90){\circle*{1.0}}
\put(110, 90){\circle*{1.0}}
\put(120, 90){\circle*{1.0}}
\put(130, 90){\circle*{1.0}}
\put(140, 90){\circle*{1.0}}
\put(150, 90){\circle*{1.0}}
\put(160, 90){\circle*{1.0}}
\put(  5,  5){\makebox(0,0){$\scriptstyle(0,0)$}}
%
\put(170,  5){\makebox(0,0){$\scriptstyle x$}}
\put(  5,110){\makebox(0,0){$\scriptstyle y$}}
\put( 20,  5){\makebox(0,0){\footnotesize\textcolor{blue}{$\scriptstyle u_4$}}}
\put( 40,  5){\makebox(0,0){\footnotesize\textcolor{blue}{$\scriptstyle u_3$}}}
\put( 60,  5){\makebox(0,0){\footnotesize\textcolor{blue}{$\scriptstyle u_2$}}}
\put( 80,  5){\makebox(0,0){\footnotesize\textcolor{blue}{$\scriptstyle u_1$}}}
\put(100,  5){\makebox(0,0){\footnotesize\textcolor{blue}{$\scriptstyle v_1$}}}
\put(120,  5){\makebox(0,0){\footnotesize\textcolor{blue}{$\scriptstyle v_2$}}}
\put(140,  5){\makebox(0,0){\footnotesize\textcolor{blue}{$\scriptstyle v_3$}}}
\put(160,  5){\makebox(0,0){\footnotesize\textcolor{blue}{$\scriptstyle v_4$}}}
\put( 25, 18){\makebox(0,0){\footnotesize\textcolor{red}{$\scriptstyle a_0$}}}
\put( 45, 18){\makebox(0,0){\footnotesize\textcolor{red}{$\scriptstyle a_0$}}}
\put( 65, 18){\makebox(0,0){\footnotesize\textcolor{red}{$\scriptstyle a_0$}}}
\put( 85, 18){\makebox(0,0){\footnotesize\textcolor{red}{$\scriptstyle a_0$}}}
\put( 95, 18){\makebox(0,0){\footnotesize\textcolor{red}{$\scriptstyle b_1$}}}
\put(115, 18){\makebox(0,0){\footnotesize\textcolor{red}{$\scriptstyle b_1$}}}
\put(135, 18){\makebox(0,0){\footnotesize\textcolor{red}{$\scriptstyle b_1$}}}
\put(155, 18){\makebox(0,0){\footnotesize\textcolor{red}{$\scriptstyle b_1$}}}
\put( 35, 28){\makebox(0,0){\footnotesize\textcolor{red}{$\scriptstyle a_1$}}}
\put( 55, 28){\makebox(0,0){\footnotesize\textcolor{red}{$\scriptstyle a_1$}}}
\put( 75, 28){\makebox(0,0){\footnotesize\textcolor{red}{$\scriptstyle a_1$}}}
\put(105, 28){\makebox(0,0){\footnotesize\textcolor{red}{$\scriptstyle b_2$}}}
\put(125, 28){\makebox(0,0){\footnotesize\textcolor{red}{$\scriptstyle b_2$}}}
\put(145, 28){\makebox(0,0){\footnotesize\textcolor{red}{$\scriptstyle b_2$}}}
\put( 45, 38){\makebox(0,0){\footnotesize\textcolor{red}{$\scriptstyle a_2$}}}
\put( 65, 38){\makebox(0,0){\footnotesize\textcolor{red}{$\scriptstyle a_2$}}}
\put( 85, 38){\makebox(0,0){\footnotesize\textcolor{red}{$\scriptstyle a_2$}}}
\put( 95, 38){\makebox(0,0){\footnotesize\textcolor{red}{$\scriptstyle b_3$}}}
\put(115, 38){\makebox(0,0){\footnotesize\textcolor{red}{$\scriptstyle b_3$}}}
\put(135, 38){\makebox(0,0){\footnotesize\textcolor{red}{$\scriptstyle b_3$}}}
\put( 55, 48){\makebox(0,0){\footnotesize\textcolor{red}{$\scriptstyle a_3$}}}
\put( 75, 48){\makebox(0,0){\footnotesize\textcolor{red}{$\scriptstyle a_3$}}}
\put(105, 48){\makebox(0,0){\footnotesize\textcolor{red}{$\scriptstyle b_4$}}}
\put(125, 48){\makebox(0,0){\footnotesize\textcolor{red}{$\scriptstyle b_4$}}}
\put( 65, 58){\makebox(0,0){\footnotesize\textcolor{red}{$\scriptstyle a_4$}}}
\put( 85, 58){\makebox(0,0){\footnotesize\textcolor{red}{$\scriptstyle a_4$}}}
\put( 95, 58){\makebox(0,0){\footnotesize\textcolor{red}{$\scriptstyle b_5$}}}
\put(115, 58){\makebox(0,0){\footnotesize\textcolor{red}{$\scriptstyle b_5$}}}
\put( 75, 68){\makebox(0,0){\footnotesize\textcolor{red}{$\scriptstyle a_5$}}}
\put(105, 68){\makebox(0,0){\footnotesize\textcolor{red}{$\scriptstyle b_6$}}}
\put( 85, 78){\makebox(0,0){\footnotesize\textcolor{red}{$\scriptstyle a_6$}}}
\put( 95, 78){\makebox(0,0){\footnotesize\textcolor{red}{$\scriptstyle b_7$}}}
\end{picture}
\caption{$t=1$ and $n=4$}\label{fig:t=1}
\end{center}
\end{figure}
The other identities can be proven similarly.
For example,
if $t=1$,
there is only one $n$-path $\boldsymbol{P}=(P_1,\dots,P_n)$
that connect $\boldsymbol{u}$ to $\boldsymbol{v}$
as in Figure~\ref{fig:t=1}.
As the product of the weights of all edges in $\boldsymbol{P}$ 
we obtain \eqref{eq:t=1}.
\input pic06.tex
If $t=2$,
there are $(n+1)$ ways to connect $\boldsymbol{u}$ to $\boldsymbol{v}$
with $n$-path $\boldsymbol{P}=(P_1,\dots,P_n)$.
As an example,
we show one way in Figure~\ref{fig:t=2}.
A similar reasoning leads to \eqref{eq:t=2}.
One can also derive \eqref{eq:t=3} by a similar argument.
\end{demo}

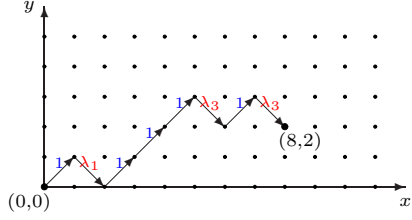
\begin{figure}[h]
\begin{center}
\setlength{\unitlength}{0.4mm}
\begin{picture}(150,80)
%
%
\put( 10, 10){\vector( 1, 0){120}}
\put( 10, 10){\vector( 0, 1){ 60}}
\put( 10, 10){\vector( 1, 1){10}}
\put( 20, 20){\vector( 1,-1){10}}
\put( 30, 10){\vector( 1, 1){10}}
\put( 40, 20){\vector( 1, 1){10}}
\put( 50, 30){\vector( 1, 1){10}}
\put( 60, 40){\vector( 1,-1){10}}
\put( 70, 30){\vector( 1, 1){10}}
\put( 80, 40){\vector( 1,-1){10}}
%
%
\put( 10, 10){\circle*{2}}
\put( 90, 30){\circle*{2}}
\put( 20, 10){\circle*{1.0}}
\put( 30, 10){\circle*{1.0}}
\put( 40, 10){\circle*{1.0}}
\put( 50, 10){\circle*{1.0}}
\put( 60, 10){\circle*{1.0}}
\put( 70, 10){\circle*{1.0}}
\put( 80, 10){\circle*{1.0}}
\put( 90, 10){\circle*{1.0}}
\put(100, 10){\circle*{1.0}}
\put(110, 10){\circle*{1.0}}
\put(120, 10){\circle*{1.0}}
\put( 10, 20){\circle*{1.0}}
\put( 20, 20){\circle*{1.0}}
\put( 30, 20){\circle*{1.0}}
\put( 40, 20){\circle*{1.0}}
\put( 50, 20){\circle*{1.0}}
\put( 60, 20){\circle*{1.0}}
\put( 70, 20){\circle*{1.0}}
\put( 80, 20){\circle*{1.0}}
\put( 90, 20){\circle*{1.0}}
\put(100, 20){\circle*{1.0}}
\put(110, 20){\circle*{1.0}}
\put(120, 20){\circle*{1.0}}
\put( 10, 30){\circle*{1.0}}
\put( 20, 30){\circle*{1.0}}
\put( 30, 30){\circle*{1.0}}
\put( 40, 30){\circle*{1.0}}
\put( 50, 30){\circle*{1.0}}
\put( 60, 30){\circle*{1.0}}
\put( 70, 30){\circle*{1.0}}
\put( 80, 30){\circle*{1.0}}
\put( 90, 30){\circle*{1.0}}
\put(100, 30){\circle*{1.0}}
\put(110, 30){\circle*{1.0}}
\put(120, 30){\circle*{1.0}}
\put( 10, 40){\circle*{1.0}}
\put( 20, 40){\circle*{1.0}}
\put( 30, 40){\circle*{1.0}}
\put( 40, 40){\circle*{1.0}}
\put( 50, 40){\circle*{1.0}}
\put( 60, 40){\circle*{1.0}}
\put( 70, 40){\circle*{1.0}}
\put( 80, 40){\circle*{1.0}}
\put( 90, 40){\circle*{1.0}}
\put(100, 40){\circle*{1.0}}
\put(110, 40){\circle*{1.0}}
\put(120, 40){\circle*{1.0}}
\put( 10, 50){\circle*{1.0}}
\put( 20, 50){\circle*{1.0}}
\put( 30, 50){\circle*{1.0}}
\put( 40, 50){\circle*{1.0}}
\put( 50, 50){\circle*{1.0}}
\put( 60, 50){\circle*{1.0}}
\put( 70, 50){\circle*{1.0}}
\put( 80, 50){\circle*{1.0}}
\put( 90, 50){\circle*{1.0}}
\put(100, 50){\circle*{1.0}}
\put(110, 50){\circle*{1.0}}
\put(120, 50){\circle*{1.0}}
\put( 10, 60){\circle*{1.0}}
\put( 20, 60){\circle*{1.0}}
\put( 30, 60){\circle*{1.0}}
\put( 40, 60){\circle*{1.0}}
\put( 50, 60){\circle*{1.0}}
\put( 60, 60){\circle*{1.0}}
\put( 70, 60){\circle*{1.0}}
\put( 80, 60){\circle*{1.0}}
\put( 90, 60){\circle*{1.0}}
\put(100, 60){\circle*{1.0}}
\put(110, 60){\circle*{1.0}}
\put(120, 60){\circle*{1.0}}
\put(  5,  5){\makebox(0,0){$\scriptstyle(0,0)$}}
\put( 95, 25){\makebox(0,0){$\scriptstyle(8,2)$}}
%
\put(130,  5){\makebox(0,0){$\scriptstyle x$}}
\put(  5, 70){\makebox(0,0){$\scriptstyle y$}}
\put( 15, 18){\makebox(0,0){\footnotesize\textcolor{blue}{$\scriptstyle1$}}}
\put( 25, 18){\makebox(0,0){\footnotesize\textcolor{red}{$\scriptstyle\lambda_{1}$}}}
\put( 35, 18){\makebox(0,0){\footnotesize\textcolor{blue}{$\scriptstyle1$}}}
\put( 45, 28){\makebox(0,0){\footnotesize\textcolor{blue}{$\scriptstyle1$}}}
\put( 55, 38){\makebox(0,0){\footnotesize\textcolor{blue}{$\scriptstyle1$}}}
\put( 65, 38){\makebox(0,0){\footnotesize\textcolor{red}{$\scriptstyle\lambda_{3}$}}}
\put( 75, 38){\makebox(0,0){\footnotesize\textcolor{blue}{$\scriptstyle1$}}}
\put( 85, 38){\makebox(0,0){\footnotesize\textcolor{red}{$\scriptstyle\lambda_{3}$}}}
%
\end{picture}
\caption{A Dyck Path of weight $\lambda_1\lambda_3^2$}\label{fig:Path2}
\end{center}
\end{figure}

\begin{figure}[hb]
\begin{center}
\setlength{\unitlength}{0.5mm}
\begin{picture}(180,60)
%
%
\put( 10, 10){\vector( 1, 0){ 60}}
\put( 10, 10){\vector( 0, 1){ 40}}
\put(100, 10){\vector( 1, 0){ 60}}
\put(100, 10){\vector( 0, 1){ 40}}
\put( 10, 10){\vector( 1, 1){10}}
\put( 20, 20){\vector( 1,-1){10}}
\put( 30, 10){\vector( 1, 1){10}}
\put( 40, 20){\vector( 1,-1){10}}
\put(100, 10){\vector( 1, 1){10}}
\put(110, 20){\vector( 1, 1){10}}
\put(120, 30){\vector( 1,-1){10}}
\put(130, 20){\vector( 1,-1){10}}
%
%
\put( 10, 10){\circle*{2}}
\put( 50, 10){\circle*{2}}
\put(100, 10){\circle*{2}}
\put(140, 10){\circle*{2}}
\put( 20, 10){\circle*{1.0}}
\put( 30, 10){\circle*{1.0}}
\put( 40, 10){\circle*{1.0}}
\put( 50, 10){\circle*{1.0}}
\put( 60, 10){\circle*{1.0}}
\put(100, 10){\circle*{1.0}}
\put(110, 10){\circle*{1.0}}
\put(120, 10){\circle*{1.0}}
\put(130, 10){\circle*{1.0}}
\put(140, 10){\circle*{1.0}}
\put(150, 10){\circle*{1.0}}
\put( 10, 20){\circle*{1.0}}
\put( 20, 20){\circle*{1.0}}
\put( 30, 20){\circle*{1.0}}
\put( 40, 20){\circle*{1.0}}
\put( 50, 20){\circle*{1.0}}
\put( 60, 20){\circle*{1.0}}
\put(100, 20){\circle*{1.0}}
\put(110, 20){\circle*{1.0}}
\put(120, 20){\circle*{1.0}}
\put(130, 20){\circle*{1.0}}
\put(140, 20){\circle*{1.0}}
\put(150, 20){\circle*{1.0}}
\put( 10, 30){\circle*{1.0}}
\put( 20, 30){\circle*{1.0}}
\put( 30, 30){\circle*{1.0}}
\put( 40, 30){\circle*{1.0}}
\put( 50, 30){\circle*{1.0}}
\put( 60, 30){\circle*{1.0}}
\put(100, 30){\circle*{1.0}}
\put(110, 30){\circle*{1.0}}
\put(120, 30){\circle*{1.0}}
\put(130, 30){\circle*{1.0}}
\put(140, 30){\circle*{1.0}}
\put(150, 30){\circle*{1.0}}
\put( 10, 40){\circle*{1.0}}
\put( 20, 40){\circle*{1.0}}
\put( 30, 40){\circle*{1.0}}
\put( 40, 40){\circle*{1.0}}
\put( 50, 40){\circle*{1.0}}
\put( 60, 40){\circle*{1.0}}
\put(100, 40){\circle*{1.0}}
\put(110, 40){\circle*{1.0}}
\put(120, 40){\circle*{1.0}}
\put(130, 40){\circle*{1.0}}
\put(140, 40){\circle*{1.0}}
\put(150, 40){\circle*{1.0}}
\put(  5,  5){\makebox(0,0){$\scriptstyle(0,0)$}}
\put( 55,  5){\makebox(0,0){$\scriptstyle(4,0)$}}
\put(105,  5){\makebox(0,0){$\scriptstyle(0,0)$}}
\put(145,  5){\makebox(0,0){$\scriptstyle(4,0)$}}
\put( 70,  5){\makebox(0,0){$\scriptstyle x$}}
\put(160,  5){\makebox(0,0){$\scriptstyle x$}}
\put(  5, 50){\makebox(0,0){$\scriptstyle y$}}
\put( 95, 50){\makebox(0,0){$\scriptstyle y$}}
\put( 15, 18){\makebox(0,0){\footnotesize\textcolor{blue}{$\scriptstyle1$}}}
\put( 25, 18){\makebox(0,0){\footnotesize\textcolor{red}{$\scriptstyle\lambda_{1}$}}}
\put( 35, 18){\makebox(0,0){\footnotesize\textcolor{blue}{$\scriptstyle1$}}}
\put( 45, 18){\makebox(0,0){\footnotesize\textcolor{red}{$\scriptstyle\lambda_{1}$}}}
\put(105, 18){\makebox(0,0){\footnotesize\textcolor{blue}{$\scriptstyle1$}}}
\put(115, 28){\makebox(0,0){\footnotesize\textcolor{blue}{$\scriptstyle1$}}}
\put(125, 28){\makebox(0,0){\footnotesize\textcolor{red}{$\scriptstyle\lambda_{2}$}}}
\put(135, 18){\makebox(0,0){\footnotesize\textcolor{red}{$\scriptstyle\lambda_{1}$}}}
\end{picture}
\caption{Dyck Paths in $\Cal{D}_{4,0}$}\label{fig:Path3}
\end{center}
\end{figure}
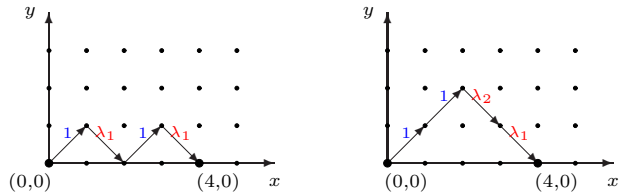

We assign the following weight to each step:
the weight of a rise vector is $1$, while
the weight of a fall vector of height $h$ is
\begin{equation}
\lambda_{h}=
\begin{cases}
\frac{q^{k}(1-aq^{k+1})(1-abq^{k+1})}{(1-abq^{2k+1})(1-abq^{2k+2})}
&\text{ if $h=2k+1$ is odd,}\\
\frac{aq^{k}(1-q^{k})(1-bq^{k})}{(1-abq^{2k})(1-abq^{2k+1})}
&\text{ if $h=2k$ is even.}
\end{cases}
\label{eq:weight}
\end{equation}
For example,
we have 
$\lambda_{1}=\frac{1-aq}{1-abq^2}$,
$\lambda_{2}=\frac{aq(1-q)(1-bq)}{(1-abq^2)(1-abq^3)}$,
$\lambda_{3}=\frac{q(1-aq^2)(1-abq^2)}{(1-abq^3)(1-abq^4)}$,
and an example of the weight of a path is Figure~\ref{fig:Path2}.
%
%
%
\begin{lemma}
Let $m$ and $n$ be a non-negative integers such that $m\equiv n$ ($\MOD\ 2$).
Then the generating function of $\Cal{D}_{m,n}$ is given by
\begin{equation}
\operatorname{GF}\left(\Cal{D}_{m,n}\right)
=\qbinom{\lfloor\frac{m}2\rfloor}{\lfloor\frac{n}2\rfloor}
\frac{(aq^{1+\lceil\frac{n}2\rceil};q)_{\frac{m-n}2}}{(abq^{2+n};q)_{\frac{m-n}2}}.
\label{eq:Dyck_GF}
\end{equation}
Here $\lfloor x\rfloor$ (resp. $\lceil x\rceil$) stands for the greatest integer that does not exceed $x$
(resp. the smallest integer that is not smaller than $x$).
Especially,
we have
\begin{equation}
\operatorname{GF}\left(\Cal{D}_{2n,0}\right)
=\frac{(aq;q)_{n}}{(abq^2;q)_{n}}.
\label{eq:Dyck_GF_special}
\end{equation}
\end{lemma}
\begin{demo}{Proof}
We prove \eqref{eq:Dyck_GF}
by induction on $m$.
If $m=0$,
then it is obvious that $\operatorname{GF}\left(\Cal{D}_{0,n}\right)$ equals $1$
if $n=0$,
and $0$ otherwise.
Assume that \eqref{eq:Dyck_GF} holds up to $m-1$.
Then we have
\begin{equation*}
\GF{\Cal{D}_{m,n}}=\GF{\Cal{D}_{m-1,n-1}}+\lambda_{n+1}\GF{\Cal{D}_{m-1,n+1}}.
\end{equation*}
If $m=2r$ and $n=2s$,
then,
by induction hyperthesis and the above recursion,
we obtain $\GF{\Cal{D}_{2r,2s}}$ equals
\begin{align*}
&\qbinom{r-1}{s-1}
\frac{\left(aq^{s+1};q\right)_{r-s}}{\left(abq^{2s+1};q\right)_{r-s}}
+\frac{q^{s}(1-aq^{s+1})(1-abq^{s+1})}{(1-abq^{2s+1})(1-abq^{2s+2})}
\qbinom{r-1}{s}
\frac{\left(aq^{s+2};q\right)_{r-s-1}}{\left(abq^{2s+3};q\right)_{r-s-1}}
\\&
=\frac{\left(q;q\right)_{r-1}\left(aq^{s+1};q\right)_{r-s}}
{\left(q;q\right)_{s}\left(q;q\right)_{r-s}\left(abq^{2s+1};q\right)_{r-s+1}}
\left\{
\left(1-q^{s}\right)\left(1-abq^{r+s+1}\right)+q^{s}\left(1-q^{r-s}\right)\left(1-abq^{s+1}\right)
\right\}
\\&
=\qbinom{r}{s}\frac{\left(aq^{s+1};q\right)_{r-s}}{\left(abq^{2s+2};q\right)_{r-s}}.
\end{align*}
%
%
\begin{figure}[ht]
\begin{center}
\setlength{\unitlength}{0.7mm}
\begin{picture}(100,50)
%
%
\put( 10, 10){\vector( 1, 0){ 60}}
\put( 10, 10){\vector( 0, 1){ 40}}
\put( 40, 20){\vector( 1, 1){ 9}}
\put( 40, 40){\vector( 1,-1){ 9}}
%
\put( 40, 20){\circle*{2}}
\put( 40, 40){\circle*{2}}
\put( 50, 30){\circle{2}}
\put( 10, 10){\circle*{1.0}}
\put( 20, 10){\circle*{1.0}}
\put( 30, 10){\circle*{1.0}}
\put( 40, 10){\circle*{1.0}}
\put( 50, 10){\circle*{1.0}}
\put( 60, 10){\circle*{1.0}}
\put( 10, 20){\circle*{1.0}}
\put( 20, 20){\circle*{1.0}}
\put( 30, 20){\circle*{1.0}}
\put( 40, 20){\circle*{1.0}}
\put( 50, 20){\circle*{1.0}}
\put( 60, 20){\circle*{1.0}}
\put( 10, 30){\circle*{1.0}}
\put( 20, 30){\circle*{1.0}}
\put( 30, 30){\circle*{1.0}}
\put( 40, 30){\circle*{1.0}}
\put( 60, 30){\circle*{1.0}}
\put( 10, 40){\circle*{1.0}}
\put( 20, 40){\circle*{1.0}}
\put( 30, 40){\circle*{1.0}}
\put( 40, 40){\circle*{1.0}}
\put( 50, 40){\circle*{1.0}}
\put( 60, 40){\circle*{1.0}}
%
%
\put( 70,  5){\makebox(0,0){$x$}}
\put(  5, 50){\makebox(0,0){$y$}}
\put( 45, 28){\makebox(0,0){\footnotesize\textcolor{blue}{$1$}}}
\put( 50, 36){\makebox(0,0){\footnotesize\textcolor{red}{$\lambda_{n+1}$}}}
\put( 59, 25){\makebox(0,0){\footnotesize\textcolor{green}{$\GF{\Cal{D}_{m,n}}$}}}
\put( 35, 15){\makebox(0,0){\footnotesize\textcolor{green}{$\GF{\Cal{D}_{m-1,n-1}}$}}}
\put( 35, 45){\makebox(0,0){\footnotesize\textcolor{green}{$\GF{\Cal{D}_{m-1,n+1}}$}}}
\end{picture}
\caption{$\GF{\Cal{D}_{m,n}}=\GF{\Cal{D}_{m-1,n-1}}+\lambda_{n+1}\GF{\Cal{D}_{m-1,n+1}}$}\label{fig:path_induction}
\end{center}
\end{figure}
This equals the right-hand side of \eqref{eq:Dyck_GF} with $m=2r$ and $n=2s$.
Hence \eqref{eq:Dyck_GF} holds when $m=2r$.
One can prove \eqref{eq:Dyck_GF} similarly when $m=2r+1$ and $n=2s+1$.
\end{demo}
For example,
if $m=4$ and $n=0$,
then $\Cal{D}_{4,0}$ has the two Dyck paths shown in Figure~\ref{fig:Path3}.
%
%
%
%
%
%
%
%
%
Thus,
the generating function of $\Cal{D}_{4,0}$ equals
\begin{align*}
\operatorname{GF}(\Cal{D}_{4,0})
=\lambda_{1}^2+\lambda_{1}\lambda_{2}
=\frac{(1-aq)(1-aq^2)}{(1-abq^2)(1-abq^3)}.
\end{align*}

\begin{demo}{Proof of Theorem~\ref{th:q-Catalan-Hankel}}
If we use
\eqref{eq:t=0}, \eqref{eq:weight} and \eqref{eq:Dyck_GF_special},
then we conclude that $\det\left(\mu_{i+j}\right)_{0\leq i,j\leq n-1}$ equals
\[
\prod_{i=1}^{n}\left(\lambda_{2i-1}\lambda_{2i}\right)^{n-i}
=\prod_{i=1}^{n}\left\{\frac{aq^{2i-1}(1-q^{i})(1-aq^{i})(1-bq^{i})(1-abq^{i})}{(1-abq^{2i-1})(1-abq^{2i})^2(1-abq^{2i+1})}\right\}^{n-i}.
\]
An easy computation leads to \eqref{eq:q-catalan-det}.
\end{demo}

\begin{remark}
One can also prove Theorem~\ref{th:q-Catalan-Hankel}
by using Motzkin paths and giving the weight $\lambda_{2h+1}$ to rise vector of hight $h$,
$\lambda_{2h}$ to fall vector of hight $h$
and $\lambda_{2h}+\lambda_{2h+1}$ to level vector of hight $h$.
Then one can prove
\begin{equation}
\operatorname{GF}\left(\Cal{M}_{m,n}\right)
=q^{\binom{m}{n}}\qbinom{m}{n}
\frac{(aq;q)_{m}(1-abq^{2n+1})}{(abq^{n+1};q)_{m+1}}.
\end{equation}
for nonnegative integers $m$ and $n$.
\end{remark}

\section{
Orthogonal Polynomials
}\label{orthogonal}

In this section we give our second proof of Theorem~\ref{th:q-Catalan-Hankel}
using the little $q$-Jacobi polynomials.
We use the notation 
$
S(t; \lambda_1,\lambda_2,\ldots)
$
for the Stieltjes-type continued fraction 
\begin{align}
\frac{1}
{\displaystyle 1-\frac{\lambda_{1}t}
{\displaystyle 1-\frac{\displaystyle\lambda_{2}t}
{ 1-\cdots
}}},
\label{eq:S-fraction}
\end{align}
and
$
J(t; b_0,b_1,b_2,\dots;\lambda_1,\lambda_2,\dots)
$
for the Jacobi-type continued fraction 
\begin{equation}
{1\over\displaystyle 1-b_0 x-
{\strut \lambda_1x^2\over\displaystyle 1-b_1x-
{\strut \lambda_2x^2\over\displaystyle
{\strut \ddots \over\displaystyle 1-b_nx-
{\strut \lambda_nx^2\over\displaystyle\ddots
}}}}}.
\label{eq:J-fraction}
\end{equation}

Given a moment sequence $\{\mu_n\}$,
we define the linear functional ${\cal L}: x^n\mapsto \mu_n$ on the vector space of
polynomials $\C[x]$. 
Then the monic polynomials $p_n(x)$ orthogonal with respect to $\cal L$ 
and of $\deg p_n(x)=n$
 satisfy a three term recurrence
relation (Favard's theorem), say
\begin{align}
p_{n+1}(x)=(x-b_n)p_{n}(x)-\lambda_n p_{n-1}(x),
\end{align}
where $p_{-1}(x)=0$ and $p_0(x)=1$.
The moment sequence $\{\mu_n\}$ is related to the coefficients $b_n$ and $\lambda_n$ by the identity:
\begin{equation}
1+\sum_{n\geq 1}\mu _nx^n=
J(t; b_0,b_1,b_2,\dots;\lambda_1,\lambda_2,\dots).
\end{equation}
Hereafter we assume $\lambda_0=\mu_0=1$ for simplicity of arguments.

Define $\Delta_n$ and $D_n(x)$ by
$$
\Delta_n=\left|
\begin{array}{cccc}
  \mu_0 & \mu_1 & \ldots & \mu_n  \\
  \mu_1 & \mu_2 & \ldots & \mu_{n+1} \\
  \vdots & \vdots & \vdots & \vdots  \\
  \mu_n & \mu_{n+1}  & \ldots & \mu_{2n}  \\
\end{array}
\right|,\qquad D_n(x)=\left|
\begin{array}{cccc}
  \mu_0 & \mu_1 & \ldots & \mu_n \\
  \mu_1 & \mu_2 & \ldots & \mu_{n+1} \\
  \vdots & \vdots & \vdots & \vdots  \\
  \mu_{n-1} & \mu_{n} & \ldots & \mu_{2n-1} \\
  1 & x  & \ldots & x^{n} \\
\end{array}
\right|.
$$
Then $p_n(x)=(\Delta_{n-1})^{-1}D_n(x)$ is the monic OPS for $\cal L$.

It is easy to see that
\begin{align}
{\cal L}(x^np_n(x))&=\frac{\Delta_n}{\Delta_{n-1}}=\lambda_n\lambda_{n-1}\ldots \lambda_1\mu_0,
\label{eq:Delta}\\
{\cal L}(x^{n+1}p_n(x))&=\frac{\chi_n}{\Delta_{n-1}}=\lambda_n\lambda_{n-1}\ldots \lambda_1\mu_0(b_0+\cdots +b_n),
\label{eq:chi}
\end{align}
where
$$
\chi_n=\left|
\begin{array}{cccc}
  \mu_0 & \mu_1 & \ldots & \mu_n \\
  \mu_1 & \mu_2 & \ldots & \mu_{n+1} \\
  \vdots & \vdots & \vdots & \vdots  \\
  \mu^{n-1} & \mu^{n} & \ldots & \mu_{2n-1} \\
  \mu_{n+1} & \mu_{n+2}  & \ldots & \mu_{2n+1} \\
\end{array}
\right|.
$$
Therefore
\begin{align}
\lambda_n=\frac{{\cal L}[p_n^2(x))]}{{\cal L}[p_{n-1}^2(x))]}
=\frac{\Delta_{n-2}\Delta_n}{\Delta_{n-1}^2},
\end{align}
and
\begin{align}
b_n=\frac{{\cal L}[xp_n^2(x))]}{{\cal L}[p_n^2(x))]}
=\frac{\chi_n}{\Delta_{n}}-\frac{\chi_{n-1}}{\Delta_{n-1}}.
\end{align}

\begin{theorem}[The Stieltjes-Rogers addition formula]
The formal power series $f(x)=\sum_{i\geq 0}a_ix^i/i!$ ($a_0=1$) has
the property that
$$
f(x+y)=\sum_{m\geq 0}\alpha_mf_m(x)f_m(y),
$$
where $\alpha_m$ is independent of $x$ and $y$ and
$$
f_m(x)={x^m\over m!}+\beta_m{x^{m+1}\over (m+1)!}+O(x^{m+2}),
$$
if and only if the formal power series $\hat f(x)=\sum_{i\geq 0}a_ix^i$
has the J-continued fraction expansion
$
J(x; b_0,b_1,b_2,\dots;\lambda_1,\lambda_2,\dots)
$
 with the parameters
$$
b_m=\beta_{m+1}-\beta_m\quad \textrm{and}\quad
 \lambda_m={\alpha_{m+1}\over \alpha_m},\quad m\geq 0.
$$
\end{theorem}
From \eqref{eq:Delta},
one can compute the Hankel determinants
\begin{equation}
\det\left(\mu_{i+j}\right)_{0\leq i,j\leq n-1}
=\Delta_{n-1}
=\mu_{0}^{n}\lambda_{1}^{n-1}\lambda_{2}^{n-2}\cdots\lambda_{n-2}^{2}\lambda_{n-1},
\end{equation}
 of \eqref{eq:q-catalan-det} by taking appropriate orthogonal polynomials $p_n(x)$.
%
%
%
%
Recall the definition of Heine's $q$-hypergeometric series
$$
{}_2\phi_1(a,b;c;q;x)=\sum_{n=0}^\infty\frac{(a;q)_n(b;q)_n}{(c;q)_n}\frac{x^n}{(q;q)_n}.
$$
The following is one of Heine's three-term contiguous relations for ${}_2\phi_1$:
$$
{}_2\phi_1(a,b;c;q;x)={}_2\phi_1(a,bq;cq;q;x)+\frac{(1-a)(c-b)}{(1-c)(1-cq)}x\,{}_2\phi_1(aq,bq;cq^2;q;x).
$$
It follows that
\begin{eqnarray*}
&&\frac{{}_2\phi_1(a,bq;cq;q;x)}{{}_2\phi_1(a,b;c;q;x)}\\
&=&S\left(x;\,
\frac{(1-a)(b-c)}{(1-c)(1-cq)},\, \frac{(1-bq)(a-cq)}{(1-cq)(1-cq^2)},\,
\frac{(1-aq)(bq-cq^2)}{(1-cq^2)(1-cq^3)},\ldots\right).
\end{eqnarray*}
Hence,
by induction,
 we can prove that
\[
\frac{{}_2\phi_1(a,bq;cq;q;x)}{{}_2\phi_1(a,b;c;q;x)}
=S(x;\lambda_1,\lambda_2,\dots),
\]
where
$$
\lambda_{2n+1}=\frac{(1-aq^n)(b-cq^n)q^n}{(1-cq^{2n})(1-cq^{2n+1})},\quad
\lambda_{2n}= \frac{(1-bq^n)(a-cq^n)q^{n-1}}{(1-cq^{2n-1})(1-cq^{2n})}.
$$
Making the substitution $b\leftarrow 1$, $a\leftarrow aq$ and $c\leftarrow abq$
into the above equation,
 we obtain
$$
\sum_{n\geq 0}\frac{(aq;q)_n}{(abq^2;q)_n}x^n=
S\left(x;\lambda_1,\lambda_2,\dots\right),
$$
where
$$
\lambda_{2n+1}= \frac{(1-aq^{n+1})(1-abq^{n+1})q^{n}}{(1-aq^{2n+1})(1-abq^{2n+2})},\qquad
\lambda_{2n}= \frac{(1-q^n)(1-bq^n)aq^{n}}{(1-abq^{2n})(1-abq^{2n+1})}.
$$
This corresponds to the little $q$-Jacobi polynomials. Indeed,
the little $q$-Jacobi polynomials
\begin{equation}
p_{n}(x;a,b;q)=\frac{(aq;q)_{n}}{(abq^{n+1};q)_{n}}(-1)^{n}q^{\binom{n}2}
{}_{2}\phi_{1}\left[
{{q^{-n}, abq^{n+1}}
\atop{aq}}
\,;\, q, xq
\right]
\end{equation}
are introduced in \cite{AA}.
The polynomials satisfy the recurrence equation
\begin{equation}
x p_{n}(x) = p_{n+1}(x) + (A_{n}+C_{n}) p_{n}(x) + A_{n-1}C_{n-1} p_{n-1}(x)
\end{equation}
where $p_{-1}(x)=0$, $p_{0}(x)=1$ and
\begin{equation}
A_{n}=\frac{q^{n}(1-aq^{n+1})(1-abq^{n+1})}
{(1-abq^{2n+1})(1-abq^{2n+2})},
\qquad
C_{n}=\frac{aq^{n}(1-q^{n})(1-bq^{n})}
{(1-abq^{2n})(1-abq^{2n+1})}.
\end{equation}
They are orthogonal with respect to the moment sequence
$\left\{\mu_{n}\right\}_{n\geq0}$ where
\begin{equation}
\mu_{n}=\frac{(aq;q)_{n}}{(abq^2;q)_{n}}.
\label{eq:orthogonal-det}
\end{equation}
For the passage from the Stieltjes-type continued fraction to
the Jacobi-type continued fraction we use the following
contraction formula:
\begin{align*}
S(x,\lambda_1,\lambda_2,\dots)
&=\frac{1}
{\displaystyle 1-\lambda_{1}t-\frac{\lambda_1\lambda_2t^2}
{\displaystyle 1-(\lambda_2+\lambda_3)t-\frac{\lambda_{3}\lambda_4 t^2}
{ 1-\cdots
}}}.
\end{align*}
Thus, by the same computation as in the former section,
we conclude that the determinant \eqref{eq:orthogonal-det} is equal to 
\eqref{eq:q-catalan-det}.
This proof gives us an insight to the determinant \eqref{eq:q-catalan-det}
from the point of view
of the classical orthogonal polynomial theory.



\section{
$q$-Dougall's formula
}\label{sec:LU-decomposition}

In this section we give our third proof of Theorem~\ref{th:q-Catalan-Hankel}
using $q$-Dougall's formula and LU-decomposition of the Hankel matrix.

First the following formula is known as $q$-Dougall's formula:
We have
\begin{align}
{}_{6}\phi_{5}\left[{{a,qa^{\frac12},-qa^{\frac12},b,c,d}
\atop{a^{\frac12},-a^{\frac12},aq/b,aq/c,aq/d}};q,\frac{aq}{bcd}\right]
=\frac{(qa,aq/bc,aq/bd,aq/cd;q)_{\infty}}
{(aq/b,aq/c,aq/d,aq/bcd;q)_{\infty}}.
\label{eq:q-Dougall}
\end{align}
provided $|aq/bcd|<1$
(see \cite[(2.7.1)]{GR}).
If we perform the substitution
 $a\leftarrow abq$, $b\leftarrow bq$, $c\leftarrow q^{-i}$ and $d\leftarrow q^{-j}$ in \eqref{eq:q-Dougall},
then we obtain
\begin{align}
{}_{6}\phi_{5}\left[
{{abq, a^{\frac12}b^{\frac12}q^{\frac32}, -a^{\frac12}b^{\frac12}q^{\frac32}, bq, q^{-i}, q^{-j}}
\atop{a^{\frac12}b^{\frac12}q^{\frac12}, -a^{\frac12}b^{\frac12}q^{\frac12}, aq, abq^{i+2}, abq^{j+2}}}
\,;\, q, aq^{i+j+1}
\right]
=\frac{\mu_{i+j}}
{\mu_{i} \mu_{j}},
\label{eq:q-hyperg-Dougall}
\end{align}
where 
$\mu_{n}
=\frac{(aq;q)_{n}}{(abq^{2};q)_{n}}$
as before.
If we use
\begin{align*}
(q^{-n};q)_{k}
=\frac{(q;q)_{n}}{(q;q)_{n-k}}
(-1)^{k}q^{\binom{k}{2}-nk}
\end{align*}
then this identity can be rewritten as
\begin{equation}
\sum_{k=0}^{\infty}
a^{k}q^{k^2}
\left[{{i}\atop{k}}\right]_{q}
\left[{{j}\atop{k}}\right]_{q}
\frac
{(q, abq, a^{\frac12}b^{\frac12}q^{\frac32}, -a^{\frac12}b^{\frac12}q^{\frac32}, bq;q)_{k}}
{(a^{\frac12}b^{\frac12}q^{\frac12}, -a^{\frac12}b^{\frac12}q^{\frac12}, aq, abq^{i+2}, abq^{j+2};q)_{k}}
=\frac{\mu_{i+j}}
{\mu_{i} \mu_{j}}.
\label{eq:LHS02}
\end{equation}
If we put
\begin{align}
&l_{ij}=\frac{\mu_{i}}{\mu_{j}}
\frac{(abq^{j+2};q)_{j}}{(abq^{i+2};q)_{j}}
\qbinom{i}{j}
\label{eq:q-L-matrix},\\
& u_{ij}=a^{i}q^{i^2}\mu_{i}\mu_{j}\,
\qbinom{j}{i}
\frac
{(q,abq,a^{\frac12}b^{\frac12}q^{\frac32},-a^{\frac12}b^{\frac12}q^{\frac32},bq;q)_{i}}
{(a^{\frac12}b^{\frac12}q^{\frac12},-a^{\frac12}b^{\frac12}q^{\frac12},aq,abq^{i+2},abq^{j+2};q)_{i}},
\label{eq:q-U-matrix}
\end{align}
then \eqref{eq:LHS02} implies
\begin{align}
\sum_{k=0}^{\infty} l_{ik} u_{kj}
=\mu_{i+j}
\label{eq:q-LU-decomp}
\end{align}
Note that
$L_{n}=\left( l_{ij} \right)_{0\leq i,j\leq n-1}$ is a lower triangular matrix
such that all main-diagonal entries are $1$,
and $U_{n}=\left( u_{ij} \right)_{0\leq i,j\leq n-1}$ is an upper-triangular
matrix with diagonal entries
\begin{align}
u_{ii}
&=a^iq^{i^2}\mu_{i}^2\,
\frac
{(q,abq,a^{\frac12}b^{\frac12}q^{\frac32},-a^{\frac12}b^{\frac12}q^{\frac32},bq;q)_{i}}
{(a^{\frac12}b^{\frac12}q^{\frac12},-a^{\frac12}b^{\frac12}q^{\frac12},aq,abq^{i+2},abq^{i+2};q)_{i}}
\nonumber\\&
=a^iq^{i^2}\,
\frac{1-abq^{2i+1}}{1-abq}\,
\frac
{(q,aq,bq,abq;q)_{i}}
{(abq^{2};q)_{2i}^2}.
\label{eq:diagonals}
\end{align}
Since $\left(\mu_{i+j}\right)_{0\leq i,j\leq n-1}=L_{n}U_{n}$,
$\det\left(\mu_{i+j}\right)_{0\leq i,j\leq n-1}$ is the product of the diagonal entries,
i.e.,
\[
\det\left(\mu_{i+j}\right)_{0\leq i,j\leq n-1}
=\prod_{i=0}^{n-1} u_{ii}.
\]
Using \eqref{eq:diagonals},
one can easily prove \eqref{eq:q-catalan-det} by a direct computation.

\begin{remark}
\label{rem:Desnanot-Jacobi}
We should note that
Corollary~\ref{cor:q-Catalan-Hankel} can be proven
by induction using the following Desnanot-Jacobi adjoint matrix theorem:
If $M$ is an $n\times n$ matrix,
then we have
\begin{equation}
\det M \det M^{1,n}_{1,n}
=\det M^{1}_{1} \det M^{n}_{n} - \det M^{1}_{n} \det M^{n}_{1},
\label{eq:Desnanot-Jacobi}
\end{equation}
where $M^{i_1,\dots,i_r}_{j_1,\dots,j_r}$ denotes the $(n-r)\times(n-r)$
submatrix obtained by removing rows $i_1,\dots,i_r$
and columns $j_1,\dots,j_r$ from $M$.

Corollary~\ref{cor:q-Catalan-Hankel} can be also proven as a special case of Theorem~\ref{th:q-kratt},
which will be proven in Section~\ref{sec:q-kratt}.
In fact,
if one puts $k_{i}=i+t$ in \eqref{eq:q-kratt},
then he obtains \eqref{eq:q-general-hankel}.
\end{remark}

\section{Miscellany}

\subsection{A proof of Theorem~\ref{th:q-kratt}}\label{sec:q-kratt}

In this subsection we give a proof of Theorem~\ref{th:q-kratt}.
Before we prove the formula,
we need to cite a lemma from \cite{K0,K}.
\begin{lemma}[Krattenthaler \cite{K0}]
Let $X_0$, $\dots$, $X_{n}$, $A_{1}$, $\dots$, $A_{n-1}$,
 and $B_{1}$, $\dots$, $B_{n-1}$ be indeterminates.
Then there holds
\begin{align}
\det\Biggl[
\prod_{l=1}^{j}\left(X_{i}+B_{l}\right)
\prod_{l=j+1}^{n-1}\left(X_{i}+A_{l}\right)
\Biggr]_{0\leq i,j\leq n-1}
= \prod_{0\leq i<j\leq n-1}
\left(X_{i} - X_{j}\right) 
\prod_{1\leq i\leq j\leq n-1}
\left(B_{i} - A_{j}\right).
\label{eq:Vanderemonde-Kratt}
\end{align}
\end{lemma}
\begin{demo}{Proof of Theorem~\ref{th:q-kratt}}
Using
\[
\mu_{n}=\frac{(aq;q)_{n}}{(abq^2;q)_{n}},
\]
we can write
\begin{align*}
\det\left(\mu_{k_{i}+j}\right)_{0\leq i,j\leq n-1}
&=\prod_{i=0}^{n-1}\frac{(aq;q)_{k_{i}}}{(abq^2;q)_{k_{i}+n-1}}
\det\left(
\prod_{l=1}^{j}\left(1-aq^{k_{i}+l}\right)
\prod_{l=j+1}^{n-1}\left(1-abq^{k_{i}+l+1}\right)
\right)_{0\leq i,j\leq n-1}
\\
&=\prod_{i=0}^{n-1}\frac{q^{(n-1)k_{i}}(aq;q)_{k_{i}}}{(abq^2;q)_{k_{i}+n-1}}
\det\left(
\prod_{l=1}^{j}\left(q^{-k_{i}}-aq^{l}\right)
\prod_{l=j+1}^{n-1}\left(q^{-k_{i}}-abq^{l+1}\right)
\right)_{0\leq i,j\leq n-1}.
\end{align*}
If we substitute $X_{i}=q^{-k_{i}}$,
$B_{l}=-aq^{l}$ and $A_{l}=-abq^{l+1}$ into \eqref{eq:Vanderemonde-Kratt},
then we see that
\begin{equation*}
\det\left(\mu_{k_{i}+j}\right)_{0\leq i,j\leq n-1}
=\prod_{i=0}^{n-1}\frac{q^{(n-1)k_{i}}(aq;q)_{k_{i}}}{(abq^2;q)_{k_{i}+n-1}}
\prod_{0\leq i<j\leq n-1}
\left(q^{-k_{i}} - q^{-k_{j}}\right) 
\prod_{1\leq i\leq j\leq n-1}
\left(abq^{j+1} - aq^{i}\right)
\end{equation*}
One can derive \eqref{eq:q-kratt} easily by a direct computation.
\end{demo}

\subsection{An addition formula for $_2F_1$}

In this subsection we give a new proof of \eqref{eq:Aigner01}
using 
an addition formula for ${}_2F_1$
and LU-decomposition of Motzkin Hankel matrices.
First,
we shall prove the following identity.
\begin{lemma}
If $i$ and $j$ are nonnegative integers,
then we have
\begin{align}\label{okinawaid}
\sum_{k\geq 0}{i\choose k}{j\choose k}
{}_2F_1\left[
{{\frac{k-i+1}{2}, \frac{k-i}{2}}
\atop{k+2}}
\,\,;\,4
\right]{}_2F_1\left[
{{\frac{k-j+1}{2}, \frac{k-j}{2}}
\atop{k+2}}
\,\,;\,4
\right]
={}_2F_1\left[
{{\frac{1-i-j}{2}, \frac{-i-j}{2}}
\atop{2}}
\,\,;\,4
\right].
\end{align}
\end{lemma}
\begin{demo}{Proof}
Recall  the quadratic transformation formula (see \cite[(3.1.5)]{GR}):
\begin{align}\label{quadratic}
(1-z)^a{}_2F_1(a,b;2b;2z)=
{}_2F_1\left(\frac{a}{2},\frac{a+1}{2};b+\frac{1}{2};\frac{z^2}{(1-z)^2}\right).
\end{align}
Applying \eqref{quadratic} with $a=k-i$, $b=k+3/2$ and $z=2$ we obtain
$$
{}_2F_1\left[
{{\frac{k-i+1}{2}, \frac{k-i}{2}}
\atop{k+2}}
\,\,;\,4
\right]=(-1)^{k-i}{}_2F_1\left[
{{k-i, k+3/2}
\atop{2k+3}}
\,\,;\,4
\right].
$$
Substituting $i$ by $j$ yields
$$
{}_2F_1\left[
{{\frac{k-j+1}{2}, \frac{k-j}{2}}
\atop{k+2}}
\,\,;\,4
\right]=(-1)^{k-j}{}_2F_1\left[
{{k-j, k+3/2}
\atop{2k+3}}
\,\,;\,4
\right].
$$
Now, applying \eqref{quadratic} with $a=-i-j$, $b=3/2$ and $z=2$ we obtain
$$
{}_2F_1\left[
{{\frac{1-i-j}{2}, \frac{-i-j}{2}}
\atop{2}}
\,\,;\,4
\right]=(-1)^{i+j}{}_2F_1\left[
{{-i-j, \frac{3}{2}}
\atop{3}}
\,\,;\,4
\right].
$$
Therefore we can rewrite \eqref{okinawaid} as follows:
\begin{align}\label{okinawaidbis}
\sum_{k\geq 0}{i\choose k}{j\choose k}
{}_2F_1\left[
{{k-i, k+3/2}
\atop{2k+3}}
\,\,;\,4
\right]
{}_2F_1\left[
{{k-j, k+3/2}
\atop{2k+3}}
\,\,;\,4
\right]
={}_2F_1\left[
{{-i-j, \frac{3}{2}}
\atop{3}}
\,\,;\,4
\right].
\end{align}
Now we recall a formula of Burchnall and Chaundy~\cite[(43)]{BC}:
\begin{align}\label{BCId}
{}_2F_1\left[
{{c-a, c-b}
\atop{c}}
\,\,;\,x
\right]&=\sum_{k\geq 0}\frac{(c-a)_k(a)_k(d)_k(c-b-d)_k}{k!(c+k-1)_k(c)_{2k}}x^{2k}\nonumber\\
&\times {}_2F_1\left[
{{c-a+k, c-b-d+k}
\atop{c+2k}}
\,\,;\,x
\right]
{}_2F_1\left[
{{c-a+k, d+k}
\atop{c+2k}}
\,\,;\,x
\right].
\end{align}
It is then easy to check that
the specialization of \eqref{BCId} with
$$
a=\frac{3}{2},\quad b=3+i+j,\quad c=3,\quad d=-j,\quad x=4
$$
yields \eqref{okinawaidbis}.
\end{demo}

\begin{demo}{Proof of \eqref{eq:Aigner01}}
Define $l_{ij}$ and $u_{ij}$ by
\begin{align}
&l_{ij}=\binom{i}{j}
{}_2F_1\left[
{{\frac{j-i+1}{2}, \frac{j-i}{2}}
\atop{j+2}}
\,\,;\,4
\right],\\
&u_{ij}=\binom{j}{i}
{}_2F_1\left[
{{\frac{i-j+1}{2}, \frac{i-j}{2}}
\atop{i+2}}
\,\,;\,4
\right].
\end{align}
Then $L_{n}=\left(l_{ij}\right)_{0\leq i,j\leq n-1}$
is a lower triangular matrix with all diagonal entries $1$,
and $U_{n}=\left(l_{ij}\right)_{0\leq i,j\leq n-1}$
is an upper triangular matrix with all diagonal entries $1$.
The formula \eqref{okinawaid} gives the LU-decomposition
of Motzkin Hankel matrix:
\[
\left(M_{ij}\right)=L_{n}U_{n}.
\]
Hence we conclude that $\det\left(M_{ij}\right)=1$.
\end{demo}

\subsection{A $q$-analogue of Schr\"{o}der numbers}

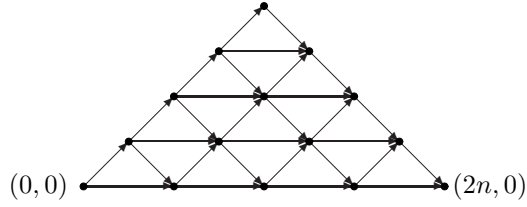
\begin{figure}[hbt]
\begin{center}
\begin{picture}(150,90)
\setlength{\unitlength}{0.6mm}
\put(0,0){\vector(1,1){10}}
\put(10,10){\vector(1,1){10}}
\put(20,20){\vector(1,1){10}}
\put(30,30){\vector(1,1){10}}
\put(20,0){\vector(1,1){10}}
\put(30,10){\vector(1,1){10}}
\put(40,20){\vector(1,1){10}}
\put(40,0){\vector(1,1){10}}
\put(50,10){\vector(1,1){10}}
\put(60,0){\vector(1,1){10}}
\put(10,10){\vector(1,-1){10}}
\put(20,20){\vector(1,-1){10}}
\put(30,30){\vector(1,-1){10}}
\put(40,40){\vector(1,-1){10}}
\put(30,10){\vector(1,-1){10}}
\put(40,20){\vector(1,-1){10}}
\put(50,30){\vector(1,-1){10}}
\put(50,10){\vector(1,-1){10}}
\put(60,20){\vector(1,-1){10}}
\put(70,10){\vector(1,-1){10}}
\put(0,0){\vector(1,0){20}}
\put(10,10){\vector(1,0){20}}
\put(20,20){\vector(1,0){20}}
\put(30,30){\vector(1,0){20}}
\put(20,0){\vector(1,0){20}}
\put(30,10){\vector(1,0){20}}
\put(40,20){\vector(1,0){20}}
\put(40,0){\vector(1,0){20}}
\put(50,10){\vector(1,0){20}}
\put(60,0){\vector(1,0){20}}
\put(0,0){\circle*{1.5}}
\put(10,10){\circle*{1.5}}
\put(20,20){\circle*{1.5}}
\put(30,30){\circle*{1.5}}
\put(40,40){\circle*{1.5}}
\put(20,0){\circle*{1.5}}
\put(30,10){\circle*{1.5}}
\put(40,20){\circle*{1.5}}
\put(50,30){\circle*{1.5}}
\put(40,0){\circle*{1.5}}
\put(50,10){\circle*{1.5}}
\put(60,20){\circle*{1.5}}
\put(60,0){\circle*{1.5}}
\put(70,10){\circle*{1.5}}
\put(80,0){\circle*{1.5}}
\put(-10,0){\makebox(0,0){$(0,0)$}}
\put(90,0){\makebox(0,0){$(2n,0)$}}
\end{picture}
\caption{Weight of Schr\"{o}der paths in $\mathcal{S}_{4,0}$}\label{fig:weight-Schroder}
\end{center}
\end{figure}

We define $S_n(q)$ ($n\geq 0$)
by the following recurrence:
\[
S_0(q)=1, \
S_n(q)=
q^{2n-1}S_{n-1}(q)+\sum_{k=0}^{n-1}q^{2(k+1)(n-1-k)}S_{n-1-k}(q)S_k(q).
\]
In fact one can show that
\[
S_n(q)=\sum_{P\in \mathcal{S}_{2n,0}}\,\omega(P),
\]
where $\omega(P)$ is the number of triangles below the path $P$ (see Figure~\ref{fig:weight-Schroder}),
and the sum runs over all Schr\"{o}der paths from the origin to $(2n,0)$.
As a $q$-analogue of \eqref{eq:Schroder,t=0,1} and \eqref{eq:Schroder,t=2}
we consider the matrix
\begin{equation}
S^{(t)}_n(q)=
\left(q^{(i-j)(i-j-1)}S_{i+j+t}(q)\right)_{0\leq i,j\leq n-1}.
\end{equation}
Note that this matrix is not a Hankel matrix,
but as a $q$-analogue of \eqref{eq:Schroder,t=0,1} and \eqref{eq:Schroder,t=2}, 
the following theorem holds:

\begin{theorem}
\label{th:q-Schroder}
Let $n$ be a positive integer.
\par\noindent
(i) If $t=0\text{ or }1$, then we have
\begin{equation}
\det S^{(1)}_{n-1}(q)=\det S^{(0)}_{n}(q)
=\prod_{k=1}^{n-1}(q^{2k-1}+1)^{n-k}.
\label{eq:Schroder01}
\end{equation}
(ii) If $t=2$,
then we have
\begin{equation}
\det S^{(2)}_n(q)=q^{-1}\prod_{k=1}^n(q^{2k-1}+1)^{n+1-k}
(\prod_{k=1}^{n+1}(q^{2k-1}+1)-1).
\label{eq:Schroder03}
\end{equation}
\end{theorem}

To prove this theorem,
we define the matrices
\begin{align*}
{\widehat S}^{(t)}_n(q)&=
\left( q^{2(n-i)(t+i+j-2)}S_{t+i+j-2}(q) \right)_{1\leq i,j\leq n},\\
{\widetilde S}^{(t)}_n(q)&=
\left( q^{-(t+i+j)(t+i+j-1)}S_{t+i+j-2}(q) \right)_{1\leq i,j\leq n},
\end{align*}
then the following lemma can be easily proven by direct computations:

\begin{lemma}
Let $n$ be a positive integer.
Then
\begin{align}
&\det {\widehat S}^{(t)}_n(q)=
q^{\frac{n(n-1)(2n+3t-4)}3}\det S_n^{(t)}(q),
\label{eqsr1}\\
&\det {\widetilde S}_n(q)=
q^{-\frac{n(4n^2+3(2t+1)n+3t^2+3t-1)}3}\det S_n^{(t)}(q).
\label{eqsr2}
\end{align}
\end{lemma}

\begin{lemma}
Let $n$ be a positive integer.
\par\noindent
(i) 
If $n\geq 2$, then we have
\begin{equation}
\det{\widehat S}^{(0)}_n(q)=q^{(n-1)(n-2)}\det{\widehat S}^{(1)}_{n-1}(q). 
\label{eqh1}
\end{equation}
(ii) 
If $n\geq 2$, then we have
\begin{equation}
\det{\widehat S}^{(1)}_n(q)=
q^{n^2-n+1}\det{\widehat S}^{(2)}_{n-1}(q)
+q^{2(n-1)^2}\det{\widehat S}^{(1)}_{n-1}(q).
\label{eqh2}
\end{equation}
(iii) 
If $n\geq 3$, then we have
\begin{equation}
\det{\widetilde S}^{(0)}_n(q)\det{\widetilde S}^{(2)}_{n-2}(q)
=
\det{\widetilde S}^{(0)}_{n-1}(q)\det{\widetilde S}^{(2)}_{n-1}(q)
-\left\{\det{\widetilde S}^{(1)}_{n-1}(q)\right\}^2.
\label{eqt1}
\end{equation}
\end{lemma}
\begin{demo}{Proof}
We consider the digraph $(V,E)$,
in which $V$ is the plane lattice $\Bbb{Z}^2$
and $E$ the set of rise vectors, fall vectors and long level vectors 
in the above half plane.
Let $u_{i}=(x_{0}-2(i-1),0)$ 
and $v^{(t)}_{j}=(x_{0}+2(j+t-1),0)$
for $i,j=1,2,\dots,n$, $t=0,1,2$ and a fixed integer $x_0$.
It is easy to see that the $n$-vertex 
$\boldsymbol{u}=(u_1,\dots,u_n)$ 
is \defterm{$D$-compatible} 
with the $n$-vertex 
$\boldsymbol{v}^{(t)}=(v^{(t)}_1,\dots,v^{(t)}_n)$.
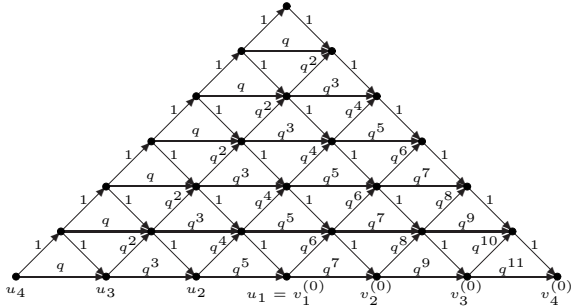
\begin{figure}[hbt]
\begin{center}
\begin{picture}(180,110)
\setlength{\unitlength}{0.6mm}
\put(0,0){\vector(1,1){10}}
\put(10,10){\vector(1,1){10}}
\put(20,20){\vector(1,1){10}}
\put(30,30){\vector(1,1){10}}
\put(40,40){\vector(1,1){10}}
\put(50,50){\vector(1,1){10}}
\put(20,0){\vector(1,1){10}}
\put(30,10){\vector(1,1){10}}
\put(40,20){\vector(1,1){10}}
\put(50,30){\vector(1,1){10}}
\put(60,40){\vector(1,1){10}}
\put(40,0){\vector(1,1){10}}
\put(50,10){\vector(1,1){10}}
\put(60,20){\vector(1,1){10}}
\put(70,30){\vector(1,1){10}}
\put(60,0){\vector(1,1){10}}
\put(70,10){\vector(1,1){10}}
\put(80,20){\vector(1,1){10}}
\put(80,0){\vector(1,1){10}}
\put(90,10){\vector(1,1){10}}
\put(100,0){\vector(1,1){10}}
\put(10,10){\vector(1,-1){10}}
\put(20,20){\vector(1,-1){10}}
\put(30,30){\vector(1,-1){10}}
\put(40,40){\vector(1,-1){10}}
\put(50,50){\vector(1,-1){10}}
\put(60,60){\vector(1,-1){10}}
\put(30,10){\vector(1,-1){10}}
\put(40,20){\vector(1,-1){10}}
\put(50,30){\vector(1,-1){10}}
\put(60,40){\vector(1,-1){10}}
\put(70,50){\vector(1,-1){10}}
\put(50,10){\vector(1,-1){10}}
\put(60,20){\vector(1,-1){10}}
\put(70,30){\vector(1,-1){10}}
\put(80,40){\vector(1,-1){10}}
\put(70,10){\vector(1,-1){10}}
\put(80,20){\vector(1,-1){10}}
\put(90,30){\vector(1,-1){10}}
\put(90,10){\vector(1,-1){10}}
\put(100,20){\vector(1,-1){10}}
\put(110,10){\vector(1,-1){10}}
\put(0,0){\vector(1,0){20}}
\put(10,10){\vector(1,0){20}}
\put(20,20){\vector(1,0){20}}
\put(30,30){\vector(1,0){20}}
\put(40,40){\vector(1,0){20}}
\put(50,50){\vector(1,0){20}}
\put(20,0){\vector(1,0){20}}
\put(30,10){\vector(1,0){20}}
\put(40,20){\vector(1,0){20}}
\put(50,30){\vector(1,0){20}}
\put(60,40){\vector(1,0){20}}
\put(40,0){\vector(1,0){20}}
\put(50,10){\vector(1,0){20}}
\put(60,20){\vector(1,0){20}}
\put(70,30){\vector(1,0){20}}
\put(60,0){\vector(1,0){20}}
\put(70,10){\vector(1,0){20}}
\put(80,20){\vector(1,0){20}}
\put(80,0){\vector(1,0){20}}
\put(90,10){\vector(1,0){20}}
\put(100,0){\vector(1,0){20}}
\put(0,0){\circle*{1.5}}
\put(10,10){\circle*{1.5}}
\put(20,20){\circle*{1.5}}
\put(30,30){\circle*{1.5}}
\put(40,40){\circle*{1.5}}
\put(50,50){\circle*{1.5}}
\put(60,60){\circle*{1.5}}
\put(20,0){\circle*{1.5}}
\put(30,10){\circle*{1.5}}
\put(40,20){\circle*{1.5}}
\put(50,30){\circle*{1.5}}
\put(60,40){\circle*{1.5}}
\put(70,50){\circle*{1.5}}
\put(40,0){\circle*{1.5}}
\put(50,10){\circle*{1.5}}
\put(60,20){\circle*{1.5}}
\put(70,30){\circle*{1.5}}
\put(80,40){\circle*{1.5}}
\put(60,0){\circle*{1.5}}
\put(70,10){\circle*{1.5}}
\put(80,20){\circle*{1.5}}
\put(90,30){\circle*{1.5}}
\put(80,0){\circle*{1.5}}
\put(90,10){\circle*{1.5}}
\put(100,20){\circle*{1.5}}
\put(100,0){\circle*{1.5}}
\put(110,10){\circle*{1.5}}
\put(120,0){\circle*{1.5}}
\put(0,-3.5){\makebox(0,0){{\tiny $u_4$}}}
\put(20,-3.5){\makebox(0,0){{\tiny $u_3$}}}
\put(40,-3.5){\makebox(0,0){{\tiny $u_2$}}}
\put(60,-3.5){\makebox(0,0){{\tiny $u_1=v^{(0)}_1$}}}
\put(80,-3.5){\makebox(0,0){{\tiny $v^{(0)}_2$}}}
\put(100,-3.5){\makebox(0,0){{\tiny $v^{(0)}_3$}}}
\put(120,-3.5){\makebox(0,0){{\tiny $v^{(0)}_4$}}}
\put(5,7){\makebox(0,0){{\tiny $1$}}}
\put(15,17){\makebox(0,0){{\tiny $1$}}}
\put(25,27){\makebox(0,0){{\tiny $1$}}}
\put(35,37){\makebox(0,0){{\tiny $1$}}}
\put(45,47){\makebox(0,0){{\tiny $1$}}}
\put(55,57){\makebox(0,0){{\tiny $1$}}}
\put(15,7){\makebox(0,0){{\tiny $1$}}}
\put(25,17){\makebox(0,0){{\tiny $1$}}}
\put(35,27){\makebox(0,0){{\tiny $1$}}}
\put(45,37){\makebox(0,0){{\tiny $1$}}}
\put(55,47){\makebox(0,0){{\tiny $1$}}}
\put(65,57){\makebox(0,0){{\tiny $1$}}}
\put(35,7){\makebox(0,0){{\tiny $1$}}}
\put(45,17){\makebox(0,0){{\tiny $1$}}}
\put(55,27){\makebox(0,0){{\tiny $1$}}}
\put(65,37){\makebox(0,0){{\tiny $1$}}}
\put(75,47){\makebox(0,0){{\tiny $1$}}}
\put(55,7){\makebox(0,0){{\tiny $1$}}}
\put(65,17){\makebox(0,0){{\tiny $1$}}}
\put(75,27){\makebox(0,0){{\tiny $1$}}}
\put(85,37){\makebox(0,0){{\tiny $1$}}}
\put(75,7){\makebox(0,0){{\tiny $1$}}}
\put(85,17){\makebox(0,0){{\tiny $1$}}}
\put(95,27){\makebox(0,0){{\tiny $1$}}}
\put(95,7){\makebox(0,0){{\tiny $1$}}}
\put(105,17){\makebox(0,0){{\tiny $1$}}}
\put(115,7){\makebox(0,0){{\tiny $1$}}}
\put(10,2){\makebox(0,0){{\tiny $q^{}$}}}
\put(20,12){\makebox(0,0){{\tiny $q^{}$}}}
\put(30,22){\makebox(0,0){{\tiny $q^{}$}}}
\put(40,32){\makebox(0,0){{\tiny $q^{}$}}}
\put(50,42){\makebox(0,0){{\tiny $q^{}$}}}
\put(60,52){\makebox(0,0){{\tiny $q^{}$}}}
\put(25,7){\makebox(0,0){{\tiny $q^{2}$}}}
\put(35,17){\makebox(0,0){{\tiny $q^{2}$}}}
\put(45,27){\makebox(0,0){{\tiny $q^{2}$}}}
\put(55,37){\makebox(0,0){{\tiny $q^{2}$}}}
\put(65,47){\makebox(0,0){{\tiny $q^{2}$}}}
\put(30,2.2){\makebox(0,0){{\tiny $q^{3}$}}}
\put(40,12.2){\makebox(0,0){{\tiny $q^{3}$}}}
\put(50,22.2){\makebox(0,0){{\tiny $q^{3}$}}}
\put(60,32.2){\makebox(0,0){{\tiny $q^{3}$}}}
\put(70,42.2){\makebox(0,0){{\tiny $q^{3}$}}}
\put(45,7){\makebox(0,0){{\tiny $q^{4}$}}}
\put(55,17){\makebox(0,0){{\tiny $q^{4}$}}}
\put(65,27){\makebox(0,0){{\tiny $q^{4}$}}}
\put(75,37){\makebox(0,0){{\tiny $q^{4}$}}}
\put(50,2.2){\makebox(0,0){{\tiny $q^{5}$}}}
\put(60,12.2){\makebox(0,0){{\tiny $q^{5}$}}}
\put(70,22.2){\makebox(0,0){{\tiny $q^{5}$}}}
\put(80,32.2){\makebox(0,0){{\tiny $q^{5}$}}}
\put(65,7){\makebox(0,0){{\tiny $q^{6}$}}}
\put(75,17){\makebox(0,0){{\tiny $q^{6}$}}}
\put(85,27){\makebox(0,0){{\tiny $q^{6}$}}}
\put(70,2.2){\makebox(0,0){{\tiny $q^{7}$}}}
\put(80,12.2){\makebox(0,0){{\tiny $q^{7}$}}}
\put(90,22.2){\makebox(0,0){{\tiny $q^{7}$}}}
\put(85,7){\makebox(0,0){{\tiny $q^{8}$}}}
\put(95,17){\makebox(0,0){{\tiny $q^{8}$}}}
\put(90,2.2){\makebox(0,0){{\tiny $q^{9}$}}}
\put(100,12.2){\makebox(0,0){{\tiny $q^{9}$}}}
\put(104,7){\makebox(0,0){{\tiny $q^{10}$}}}
\put(110,2.2){\makebox(0,0){{\tiny $q^{11}$}}}
\end{picture}
\caption{Weight of each edge}\label{fig:weight-edge}
\end{center}
\end{figure}
We assign the weight of each edge as
a rise vector, a fall vector and a long level vector 
whose origin is $(x,y)$ and ends at $(x+1,y+1)$, $(x+1,y-1)$ and $(x+2,y)$
has weight $q^{x-y-x_{0}+2(n-1)}$, $1$ and 
$q^{x-y-x_{0}+2n-1}$, respectively, 
which is visualized in Figure~\ref{fig:weight-edge}. 
Then, by applying Lemma 2.3, we can obtain
\begin{equation}
\GF{\mathrsfs{P}_0\left(\boldsymbol{u},\boldsymbol{v}^{(t)}\right)}
=\det\widehat{S}^{(t)}_n(q).
\label{eq:eq_GFofHatS}
\end{equation}
This is important to prove the following.
\par\noindent
(i) 
Assume $t=0$ and let $\boldsymbol{u}$ and $\boldsymbol{v}$ be as above.
Put  ${\widetilde u}_{i}=(x_{0}-2i+3,1)$ 
and ${\widetilde v}^{(1)}_{j}=(x_{0}+2j-3,1)$
for $i,j=2,\dots,n$,
and let $\widetilde{\boldsymbol{u}}=(\widetilde{u}_{2},\dots,\widetilde{u}_{n})$
and $\widetilde{\boldsymbol{v}}^{(1)}=(\widetilde{v}^{(1)}_{2},\dots,\widetilde{v}^{(1)}_{n})$.
Then each $n$-path $\boldsymbol{P}=(P_1,P_2,\dots,P_n)$ 
from $\boldsymbol{u}$ to $\boldsymbol{v}^{(0)}$
corresponds to an $(n-1)$-path $\widetilde{\boldsymbol{P}}=(\widetilde{P}_2,\dots,\widetilde{P}_n)$
from $\widetilde{\boldsymbol{u}}$ to $\widetilde{\boldsymbol{v}}^{(1)}$
by regarding $\widetilde{\boldsymbol{P}}$ as the subpath of $\boldsymbol{P}$.
In fact, note that $P_{1}$ is always the path composed of a single vertex $u_{1}=v_{1}$,
each $P_{i}$ always starts from the rise vector $u_{i}\rightarrow\widetilde{u}_{i}$
and ends at the fall vector $\widetilde{v}^{(1)}_{i}\rightarrow v_i^{(0)}$
for $i=2,\dots,n$.
Hence this gives a bijection,
and the product of the weight of 
the rise vectors $u_{i}\rightarrow\widetilde{u}_{i}$ and 
the fall vectors $\widetilde{v}^{(1)}_{i}\rightarrow v_i^{(0)}$ 
for $i=2,\dots,n$ is $q^{(n-1)(n-2)}$.
This proves \eqref{eqh1}.
\par\noindent
(ii) 
Assume $t=1$ and let $\boldsymbol{u}$ and $\boldsymbol{v}$ be as above,
i.e.,
$u_{i}=(x_{0}-2(i-1),0)$ 
and $v^{(1)}_{j}=(x_{0}+2j,0)$ for $1\leq i,j\leq n$
(see Figure~\ref{fig:lem-proof}).
\begin{figure}[hbt]
\begin{center}
\begin{picture}(200,130)
\setlength{\unitlength}{0.6mm}
\put(10,10){\vector(1,1){10}}
\put(20,20){\vector(1,1){10}}
\put(30,30){\vector(1,1){10}}
\put(40,40){\vector(1,1){10}}
\put(50,50){\vector(1,1){10}}
\put(60,60){\vector(1,1){10}}
\put(30,10){\vector(1,1){10}}
\put(40,20){\vector(1,1){10}}
\put(50,30){\vector(1,1){10}}
\put(60,40){\vector(1,1){10}}
\put(70,50){\vector(1,1){10}}
\put(50,10){\vector(1,1){10}}
\put(60,20){\vector(1,1){10}}
\put(70,30){\vector(1,1){10}}
\put(80,40){\vector(1,1){10}}
\put(60,0){\vector(1,1){10}}
\put(70,10){\vector(1,1){10}}
\put(80,20){\vector(1,1){10}}
\put(90,30){\vector(1,1){10}}
\put(80,0){\vector(1,1){10}}
\put(90,10){\vector(1,1){10}}
\put(100,20){\vector(1,1){10}}
\put(100,0){\vector(1,1){10}}
\put(110,10){\vector(1,1){10}}
\put(120,0){\vector(1,1){10}}
\put(10,10){\vector(1,-1){10}}
\put(20,20){\vector(1,-1){10}}
\put(30,30){\vector(1,-1){10}}
\put(40,40){\vector(1,-1){10}}
\put(50,50){\vector(1,-1){10}}
\put(60,60){\vector(1,-1){10}}
\put(70,70){\vector(1,-1){10}}
\put(30,10){\vector(1,-1){10}}
\put(40,20){\vector(1,-1){10}}
\put(50,30){\vector(1,-1){10}}
\put(60,40){\vector(1,-1){10}}
\put(70,50){\vector(1,-1){10}}
\put(80,60){\vector(1,-1){10}}
\put(50,10){\vector(1,-1){10}}
\put(60,20){\vector(1,-1){10}}
\put(70,30){\vector(1,-1){10}}
\put(80,40){\vector(1,-1){10}}
\put(90,50){\vector(1,-1){10}}
\put(70,10){\vector(1,-1){10}}
\put(80,20){\vector(1,-1){10}}
\put(90,30){\vector(1,-1){10}}
\put(100,40){\vector(1,-1){10}}
\put(100,20){\vector(1,-1){10}}
\put(110,30){\vector(1,-1){10}}
\put(120,20){\vector(1,-1){10}}
\put(0,0){\vector(1,0){20}}
\put(10,10){\vector(1,0){20}}
\put(20,20){\vector(1,0){20}}
\put(30,30){\vector(1,0){20}}
\put(40,40){\vector(1,0){20}}
\put(50,50){\vector(1,0){20}}
\put(60,60){\vector(1,0){20}}
\put(20,0){\vector(1,0){20}}
\put(30,10){\vector(1,0){20}}
\put(40,20){\vector(1,0){20}}
\put(50,30){\vector(1,0){20}}
\put(60,40){\vector(1,0){20}}
\put(70,50){\vector(1,0){20}}
\put(40,0){\vector(1,0){20}}
\put(50,10){\vector(1,0){20}}
\put(60,20){\vector(1,0){20}}
\put(70,30){\vector(1,0){20}}
\put(80,40){\vector(1,0){20}}
\put(70,10){\vector(1,0){20}}
\put(80,20){\vector(1,0){20}}
\put(90,30){\vector(1,0){20}}
\put(80,0){\vector(1,0){20}}
\put(90,10){\vector(1,0){20}}
\put(100,20){\vector(1,0){20}}
\put(100,0){\vector(1,0){20}}
\put(110,10){\vector(1,0){20}}
\put(120,0){\vector(1,0){20}}
\put(0,0){\vector(1,1){10}}
\put(20,0){\vector(1,1){10}}
\put(40,0){\vector(1,1){10}}
\put(90,10){\vector(1,-1){10}}
\put(110,10){\vector(1,-1){10}}
\put(130,10){\vector(1,-1){10}}
\put(60,0){\vector(1,0){20}}
%
\put(0,0){\circle*{1.5}}
\put(10,10){\circle*{1.5}}
\put(20,20){\circle*{1.5}}
\put(30,30){\circle*{1.5}}
\put(40,40){\circle*{1.5}}
\put(50,50){\circle*{1.5}}
\put(60,60){\circle*{1.5}}
\put(70,70){\circle*{1.5}}
\put(20,0){\circle*{1.5}}
\put(30,10){\circle*{1.5}}
\put(40,20){\circle*{1.5}}
\put(50,30){\circle*{1.5}}
\put(60,40){\circle*{1.5}}
\put(70,50){\circle*{1.5}}
\put(80,60){\circle*{1.5}}
\put(40,0){\circle*{1.5}}
\put(50,10){\circle*{1.5}}
\put(60,20){\circle*{1.5}}
\put(70,30){\circle*{1.5}}
\put(80,40){\circle*{1.5}}
\put(90,50){\circle*{1.5}}
\put(60,0){\circle*{1.5}}
\put(70,10){\circle*{1.5}}
\put(80,20){\circle*{1.5}}
\put(90,30){\circle*{1.5}}
\put(100,40){\circle*{1.5}}
\put(80,0){\circle*{1.5}}
\put(90,10){\circle*{1.5}}
\put(100,20){\circle*{1.5}}
\put(110,30){\circle*{1.5}}
\put(100,0){\circle*{1.5}}
\put(110,10){\circle*{1.5}}
\put(120,20){\circle*{1.5}}
\put(120,0){\circle*{1.5}}
\put(130,10){\circle*{1.5}}
\put(140,0){\circle*{1.5}}
\put(0,-3.5){\makebox(0,0){{\tiny $u_4$}}}
\put(20,-3.5){\makebox(0,0){{\tiny $u_3$}}}
\put(40,-3.5){\makebox(0,0){{\tiny $u_2$}}}
\put(60,-3.5){\makebox(0,0){{\tiny $u_1$}}}
\put(80,-3.5){\makebox(0,0){{\tiny $v^{(1)}_1$}}}
\put(100,-3.5){\makebox(0,0){{\tiny $v^{(1)}_2$}}}
\put(120,-3.5){\makebox(0,0){{\tiny $v^{(1)}_3$}}}
\put(140,-3.5){\makebox(0,0){{\tiny $v^{(1)}_4$}}}
\put(5,7){\makebox(0,0){{\tiny $1$}}}
\put(15,17){\makebox(0,0){{\tiny $1$}}}
\put(25,27){\makebox(0,0){{\tiny $1$}}}
\put(35,37){\makebox(0,0){{\tiny $1$}}}
\put(45,47){\makebox(0,0){{\tiny $1$}}}
\put(55,57){\makebox(0,0){{\tiny $1$}}}
\put(65,67){\makebox(0,0){{\tiny $1$}}}
\put(15,7){\makebox(0,0){{\tiny $1$}}}
\put(25,17){\makebox(0,0){{\tiny $1$}}}
\put(35,27){\makebox(0,0){{\tiny $1$}}}
\put(45,37){\makebox(0,0){{\tiny $1$}}}
\put(55,47){\makebox(0,0){{\tiny $1$}}}
\put(65,57){\makebox(0,0){{\tiny $1$}}}
\put(75,67){\makebox(0,0){{\tiny $1$}}}
\put(35,7){\makebox(0,0){{\tiny $1$}}}
\put(45,17){\makebox(0,0){{\tiny $1$}}}
\put(55,27){\makebox(0,0){{\tiny $1$}}}
\put(65,37){\makebox(0,0){{\tiny $1$}}}
\put(75,47){\makebox(0,0){{\tiny $1$}}}
\put(85,57){\makebox(0,0){{\tiny $1$}}}
\put(55,7){\makebox(0,0){{\tiny $1$}}}
\put(65,17){\makebox(0,0){{\tiny $1$}}}
\put(75,27){\makebox(0,0){{\tiny $1$}}}
\put(85,37){\makebox(0,0){{\tiny $1$}}}
\put(95,47){\makebox(0,0){{\tiny $1$}}}
\put(75,7){\makebox(0,0){{\tiny $1$}}}
\put(85,17){\makebox(0,0){{\tiny $1$}}}
\put(95,27){\makebox(0,0){{\tiny $1$}}}
\put(105,37){\makebox(0,0){{\tiny $1$}}}
\put(95,7){\makebox(0,0){{\tiny $1$}}}
\put(105,17){\makebox(0,0){{\tiny $1$}}}
\put(115,27){\makebox(0,0){{\tiny $1$}}}
\put(115,7){\makebox(0,0){{\tiny $1$}}}
\put(125,17){\makebox(0,0){{\tiny $1$}}}
\put(135,7){\makebox(0,0){{\tiny $1$}}}
\put(10,2){\makebox(0,0){{\tiny $q^{}$}}}
\put(20,12){\makebox(0,0){{\tiny $q^{}$}}}
\put(30,22){\makebox(0,0){{\tiny $q^{}$}}}
\put(40,32){\makebox(0,0){{\tiny $q^{}$}}}
\put(50,42){\makebox(0,0){{\tiny $q^{}$}}}
\put(60,52){\makebox(0,0){{\tiny $q^{}$}}}
\put(70,62){\makebox(0,0){{\tiny $q^{}$}}}
\put(25,7){\makebox(0,0){{\tiny $q^{2}$}}}
\put(35,17){\makebox(0,0){{\tiny $q^{2}$}}}
\put(45,27){\makebox(0,0){{\tiny $q^{2}$}}}
\put(55,37){\makebox(0,0){{\tiny $q^{2}$}}}
\put(65,47){\makebox(0,0){{\tiny $q^{2}$}}}
\put(75,57){\makebox(0,0){{\tiny $q^{2}$}}}
\put(30,2.2){\makebox(0,0){{\tiny $q^{3}$}}}
\put(40,12.2){\makebox(0,0){{\tiny $q^{3}$}}}
\put(50,22.2){\makebox(0,0){{\tiny $q^{3}$}}}
\put(60,32.2){\makebox(0,0){{\tiny $q^{3}$}}}
\put(70,42.2){\makebox(0,0){{\tiny $q^{3}$}}}
\put(80,52.2){\makebox(0,0){{\tiny $q^{3}$}}}
\put(45,7){\makebox(0,0){{\tiny $q^{4}$}}}
\put(55,17){\makebox(0,0){{\tiny $q^{4}$}}}
\put(65,27){\makebox(0,0){{\tiny $q^{4}$}}}
\put(75,37){\makebox(0,0){{\tiny $q^{4}$}}}
\put(85,47){\makebox(0,0){{\tiny $q^{4}$}}}
\put(50,2.2){\makebox(0,0){{\tiny $q^{5}$}}}
\put(60,12.2){\makebox(0,0){{\tiny $q^{5}$}}}
\put(70,22.2){\makebox(0,0){{\tiny $q^{5}$}}}
\put(80,32.2){\makebox(0,0){{\tiny $q^{5}$}}}
\put(90,42.2){\makebox(0,0){{\tiny $q^{5}$}}}
\put(65,7){\makebox(0,0){{\tiny $q^{6}$}}}
\put(75,17){\makebox(0,0){{\tiny $q^{6}$}}}
\put(85,27){\makebox(0,0){{\tiny $q^{6}$}}}
\put(95,37){\makebox(0,0){{\tiny $q^{6}$}}}
\put(70,2.2){\makebox(0,0){{\tiny $q^{7}$}}}
\put(80,12.2){\makebox(0,0){{\tiny $q^{7}$}}}
\put(90,22.2){\makebox(0,0){{\tiny $q^{7}$}}}
\put(100,32.2){\makebox(0,0){{\tiny $q^{7}$}}}
\put(85,7){\makebox(0,0){{\tiny $q^{8}$}}}
\put(95,17){\makebox(0,0){{\tiny $q^{8}$}}}
\put(105,27){\makebox(0,0){{\tiny $q^{8}$}}}
\put(90,2.2){\makebox(0,0){{\tiny $q^{9}$}}}
\put(100,12.2){\makebox(0,0){{\tiny $q^{9}$}}}
\put(110,22.2){\makebox(0,0){{\tiny $q^{9}$}}}
\put(104,7){\makebox(0,0){{\tiny $q^{10}$}}}
\put(114,17){\makebox(0,0){{\tiny $q^{10}$}}}
\put(110,2.2){\makebox(0,0){{\tiny $q^{11}$}}}
\put(120,12.2){\makebox(0,0){{\tiny $q^{11}$}}}
\put(124,7){\makebox(0,0){{\tiny $q^{12}$}}}
\put(130,2.2){\makebox(0,0){{\tiny $q^{13}$}}}
\end{picture}
\caption{$t=1$ case}\label{fig:lem-proof}
\end{center}
\end{figure}
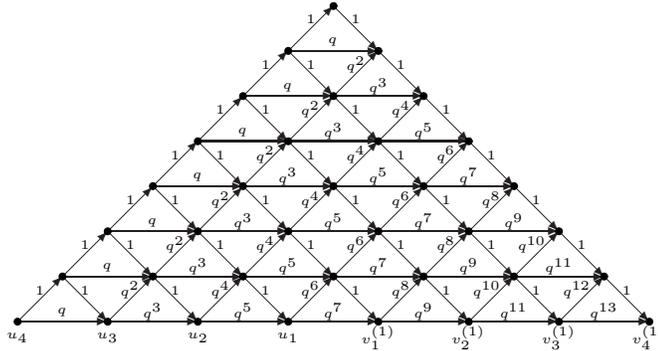
Put  ${\widetilde u}_{i}=(x_{0}-2i+3,1)$ ($2\leq i\leq n$)
and ${\widetilde v}^{(2)}_{j}=(x_{0}+2j-1,1)$
 ($2\leq j\leq n$),
and let $\widetilde{\boldsymbol{u}}=(\widetilde{u}_{2},\dots,\widetilde{u}_{n})$
and $\widetilde{\boldsymbol{v}}^{(2)}=(\widetilde{v}^{(2)}_{2},\dots,\widetilde{v}^{(2)}_{n})$.
Further,
 put ${\check u}_{i}=(x_{0}-2i+4,2)$ ($2\leq i\leq n$)
and ${\check v}^{(1)}_{j}=(x_{0}+2j-2,2)$ ($2\leq j\leq n$),
and let $\check{\boldsymbol{u}}=(\check{u}_{2},\dots,\check{u}_{n})$
and $\check{\boldsymbol{v}}^{(1)}=(\check{v}^{(1)}_{2},\dots,\check{v}^{(1)}_{n})$.
Let $\boldsymbol{P}=(P_1,P_2,\dots,P_n)$ be 
any non-intersecting $n$-paths from 
$\boldsymbol{u}$ to $\boldsymbol{v}^{(1)}$. 
Then, it is easy to see that $\boldsymbol{P}$ 
must satisfy one of the following two conditions:
\begin{enumerate}
\item[(1)] 
$P_1$ is the long level vector whose origin is $u_1$ and ends at 
$v_1^{(1)}$, 
and $P_i$ goes through the vertices 
${\widetilde u}_i$ and ${\widetilde v}_i^{(2)}$ for $i=2,3,\dots,n$. 
\item[(2)] 
$P_1$ is a path which goes through only three vertices 
$u_1$, $u_1+(1,1)(=v_1^{(1)}-(1,-1))$ and $v_1^{(1)}$, 
and $P_i$ goes through the vertices 
$\widetilde{u}_i$, $\check{u}_i$, $\check{v}_i^{(1)}$ and $\widetilde{v}_i^{(2)}$ 
for $i=2,3,\dots,n$. 
\end{enumerate}
By a similar argument as in the proof of (i), we can deduce that
\[
\GF{\mathrsfs{P}_0(\boldsymbol{u}_n,\boldsymbol{v}^{(1)}_n)}
=
q^{n^2-n+1}
\GF{\mathrsfs{P}_0(\boldsymbol{u}_{n-1},
\boldsymbol{v}^{(2)}_{n-1})}
+q^{2(n-1)^2}
\GF{\mathrsfs{P}_0(\boldsymbol{u}_{n-1},
\boldsymbol{v}^{(1)}_{n-1})}
\]
holds.
By the equality \eqref{eq:eq_GFofHatS}, we obtain the identity 
\eqref{eqh2}. 
\par\noindent
(iii)
This identity can be proven 
by applying the Desnanot-Jacobi adjoint matrix theorem 
\eqref{eq:Desnanot-Jacobi} to ${\widetilde S}^{(t)}_{n}(q)$.
\end{demo}
\begin{demo}{Proof of Theorem~\ref{th:q-Schroder}}
(i) The first equality of \eqref{eq:Schroder01} is easily obtained from 
\eqref{eqsr1} and \eqref{eqh1}. 
By applying the equalities \eqref{eqsr1} and \eqref{eqsr2} 
to \eqref{eqh2} and \eqref{eqt1}, we have 
\begin{align}
\det S^{(1)}_n(q)&=q\det S^{(2)}_{n-1}(q)+\det S^{(1)}_{n-1}(q)
\label{eqesr1}
\end{align}
for $n\geq 2$, and we have
\begin{align}
\det S^{(0)}_n(q)\det S^{(2)}_{n-2}(q)
&=
\det S^{(0)}_{n-1}(q)\det S^{(2)}_{n-1}(q)
-q^{2(n-1)}\{\det S^{(1)}_{n-1}(q)\}^2
\label{eqesr2}
\end{align}
for $n\geq 3$. 
By the equalities \eqref{eqesr1} and \eqref{eqesr2}, 
for $n\geq 3$, the following identity holds: 
\begin{align}
&\det S^{(0)}_n(q)(\det S^{(1)}_{n-1}(q)-\det S^{(1)}_{n-2}(q))
\nonumber\\
&=
\det S^{(0)}_{n-1}(q)(\det S^{(1)}_n(q)-\det S^{(1)}_{n-1}(q))
-q^{2n-1}\{\det S^{(1)}_{n-1}(q)\}^2.
\label{eqabc}
\end{align}
Moreover, by applying the first equality of \eqref{eq:Schroder01} to 
\eqref{eqabc} and replacing $n$ with $n-1$, we obtain
\begin{equation}
(1+q^{2n-3})\{\det S^{(0)}_{n-1}(q)\}^2
=
\det S^{(0)}_{n-2}(q)\det S^{(0)}_{n}(q)
\label{eqsrkey}
\end{equation}
for $n\geq 4$. 
We prove the second equality of \eqref{eq:Schroder01} by induction on $n$. 
If $n=1,2,3$, then it is easily obtained by direct computations. 
Assume that \eqref{eq:Schroder01} holds up to $n-1$. Then, 
by \eqref{eqsrkey} and induction hypothesis, we can 
obtain the second equality \eqref{eq:Schroder01}. 
\par\noindent
(ii) It follows from our result of (i) and the equality \eqref{eqesr1} 
that \eqref{eq:Schroder03} holds.
\end{demo}
%
%
%

By applying the Desnanot-Jacobi adjoint matrix theorem 
\eqref{eq:Desnanot-Jacobi} to $S^{(1)}_{n+1}(1)$, then 
we have 
\begin{equation}
\det S^{(1)}_{n+1}(1)
\det S^{(3)}_{n-1}(1)=
\det S^{(1)}_{n}(1)
\det S^{(3)}_{n}(1)-
\{\det S^{(2)}_{n}(1)\}^2
\label{eqDJ}
\end{equation}
for $n\geq 2$. 
Therefore the following 
identity is easily obtained by induction on $n$ and the formula 
\eqref{eqDJ}:
\begin{remark}
For positive integer $n$, we have
\begin{equation}
\det S^{(3)}_n(1)=
2^{\binom{n+3}{2}}-(2n+3)2^{\binom{n+2}{2}}-2^{\binom{n+1}{2}}.
\end{equation}
\end{remark}
Note that there is a relation between domino tilings of the Aztec diamonds
and Schr\"{o}der paths
(see \cite{EF,Stan}).
It might be an interesting problem to consider what this weight means.
%
%
%
%
\subsection{Delannoy numbers}

The \newterm{Delannoy numbers} $D(a,b)$ are the number of lattice paths from $(0,0)$ to $(b,a)$ 
in which only east $(1, 0)$, north $(0, 1)$, and northeast $(1, 1)$ steps are allowed. 
They are given by the recurrence relation
\begin{align}
&D(a,0)=D(0,b)=1, \nonumber\\
&D(a,b)=D(a-1,b)+D(a-1,b-1)+D(a,b-1).
\end{align}
The first few terms of $D(n,n)$ ($n=0,1,2,\dots$) are given by
$1$, $3$, $13$, $63$, $321$, $\dots$.
By a similar argument we can derive the following result.
We may give a proof in another occasion.
\begin{prop}
Let $n$ be a positive integers.
Then the following identities would hold:
\begin{align}
&\det\left(D(i+j,i+j)\right)_{0\leq i,j\leq n-1}=2^{\binom{n+1}{2}-1},\\
&\det\left(D(i+j+1,i+j+1)\right)_{0\leq i,j\leq n-1}
=2^{\binom{n+2}{2}-2}+2^{\binom{n+1}{2}-1},\\
&\det\left(D(i+j+2,i+j+2)\right)_{0\leq i,j\leq n-1}
=2^{\binom{n+3}{2}-3}+(2n+1)2^{\binom{n+2}{2}-2}-2^{\binom{n+1}{2}-1}.
\end{align}
\end{prop}

\section{Concluding remarks}

Since a hyperpfaffian version of \eqref{eq:Catalan-det} is obtained
in \cite{LT},
we believe it will be interesting problem to consider a hyperpfaffian version of 
Theorem~\ref{th:q-Catalan-Hankel} and Theorem~\ref{th:q-kratt}.
We shall argue on it in another chance.
The authors also would like to express their gratitude
to the anonymous referee for his (her) constructive comments.


\end{document}